\documentclass{amsart}
\usepackage{amsfonts,amssymb,amsthm,amsmath,enumerate}

\usepackage{xypic}

\theoremstyle{plain}

\newtheorem*{Th*}{Theorem}

\newtheorem*{Cor*}{Corollary}

\theoremstyle{definition}

\theoremstyle{remark}

\numberwithin{equation}{section}

\makeatletter
\def\Set@Scallop[#1]#2#3{{#1}\Parens{#2}{#3}}
\newcommand\DeclareScalableOperator[2]{%
  \expandafter\def\csname#1\endcsname{\@ifnextchar[{{#2}\Set@Scallop}{{#2}\Set@Scallop[{}]}}
}
\makeatother

\newcommand{\fa}{for all }

\newcommand{\fs}{for some }
\newcommand\mathfa[1][{}]{\quad\text{\fa{#1} }}

\newcommand{\scth}{such that }
\newcommand{\AND}{and}

\newcommand\mathtxt[1]{\quad\text{{#1}}\quad}

\newcommand{\nd}{\mathtxt\AND}

\newcommand\cf{\emph{cf.}~}
\newcommand\vq{\emph{v.}~}

\newcommand\ie{\emph{i.e.}~}
\newcommand\viz{\emph{viz.}~}
\newcommand\via{\emph{via}~}
\newcommand\loccit{\emph{loc.~cit.}}
\newcommand\opcit{\emph{op.~cit.}}

\newcommand\etc{\emph{etc.}}

\newcommand\vphi{\varphi}

\newcommand\eps{\varepsilon}
\newcommand\nats{\mathbb{N}}
\newcommand\ints{\mathbb{Z}}

\newcommand\reals{\mathbb{R}}
\newcommand\cplxs{\mathbb{C}}
\newcommand\knums{\mathbb K}

\newcommand\vvoid{\varnothing}
\newcommand\sle{\leqslant}
\newcommand\sge{\geqslant}

\DeclareMathOperator\Ad{\mathrm{Ad}}

\DeclareMathOperator\id{\mathrm{id}}

\DeclareMathOperator\ind{\mathrm{ind}}

\makeatletter
\newcommand\Size[7][1]{
                                 \ifx#20%
                                        \def\r@l{}\def\r@m{}\def\r@r{}%
                                 \else%
                                    \ifx#21%
                                           \def\r@l{\bigl}\def\r@r{\bigr}\def\r@m{\bigm}%
                                    \else%
                                           \ifx#22%
                                                 \def\r@l{\Bigl}\def\r@r{\Bigr}\def\r@m{\Bigm}%
                                            \else%
                                                 \ifx#23%
                                                        \def\r@l{\biggl}\def\r@r{\biggr}\def\r@m{\biggm}%
                                                  \else
                                                        \ifx#24%
                                                        \def\r@l{\Biggl}\def\r@r{\Biggr}\def\r@m{\Biggm}%
                                                        \fi%
                                                  \fi%
                                            \fi%
                                      \fi%
                                 \fi%
                                 \ifx#10%
                                       \def\r@m{}%
                                 \fi%
                                 \r@l#3{#4}\r@m#5{#6}\r@r#7%
}%

\makeatother

\newcommand\Set[3]{
                                 \Size{#1}{\{}{#2}{|}{#3}{\}}%
}%
\newcommand\Dual[3]{
                                 \Size[0]{#1}{\langle}{#2}{,}{#3}{\rangle}%
}%
\newcommand\Parens[2]{
  \Size[0]{#1}{(}{#2}{}{}{)}
}

\newcommand\Abs[2]{
  \Size[0]{#1}{\lvert}{#2}{}{}{\rvert}
}

\makeatletter

\newcommand{\IfUpperCase}[1]{\begingroup 
  \protected@edef\@tempa{\expandafter\@firstofone\@firstofone#1.}%
  \expandafter\IfUpperCasE \@tempa\delimiter}

\def\IfUpperCasE #1#2\delimiter{%
  \protected@edef\@tempa{\meaning#1\meaning a}%
  \ifnum \expandafter\IfUppercaSE\@tempa \IfUppercaSE
   \endgroup \expandafter\@firstoftwo
  \else
   \endgroup \expandafter\@secondoftwo
  \fi}

\@ifundefined{strip@prefix}{\def\strip@prefix#1>{}}{}
\def\@tempa{the letter }
\edef\@tempa{\expandafter\strip@prefix\meaning\@tempa}

\expandafter\def\expandafter\IfUppercaSE\expandafter#\expandafter1\@tempa#2#3\IfUppercaSE{\uccode`#2=`#2 }

\newif\ifuc@se
\def\setuc@se#1{\IfUpperCase{#1}{\uc@setrue}{\uc@sefalse}}

\def\theoremn@me#1{\ifuc@se \lowercase{\csname#1name\endcsname}\ignorespaces%
  \else \edef\@temp{\lowercase{\lowercase{\csname#1name\endcsname}}}\@temp\ignorespaces%
  \fi}
\def\theoremn@mes#1{\ifuc@se \lowercase{\csname#1names\endcsname}\ignorespaces%
  \else \edef\@temp{\lowercase{\lowercase{\csname#1names\endcsname}}}\@temp\ignorespaces%
  \fi}
  
\def\thmref#1#2{\setuc@se{#1}\lowercase{{\theoremn@me{#1}\lowercase{\ref{#1:#2}}}}}

\newcommand{\DefTheorem}[2]{\newenvironment{#1}[1][\empty]{\ignorespaces\begin{#2}\ifx##1\empty{}\else\lowercase{\label{#1:##1}}\fi\ignorespaces}{\end{#2}\ignorespacesafterend}}

\DefTheorem{Th}{theorem}
\DefTheorem{Prop}{proposition}
\DefTheorem{Cor}{corollary}
\DefTheorem{Lem}{lemma}
\DefTheorem{Def}{definition}
\DefTheorem{Rem}{remark}
\DefTheorem{Par}{para}

\newenvironment{Par*}{\ignorespaces\noindent\ignorespaces}{\ignorespacesafterend}

\makeatletter

\newif\if@smallmat
\newif\if@none
\newif\if@paren
\newif\if@brack
\newif\if@brace
\newif\if@vline
                                 {\ifx#20%
                                        \@smallmattrue%
                                  \else%
                                         \@smallmatfalse
                                  \fi%
                                  \ifx#11%
                                         \@nonefalse\@parentrue\@brackfalse\@bracefalse\@vlinefalse%
                                  \else%
                                       \ifx#12%
                                            \@nonefalse\@parenfalse\@bracktrue\@bracefalse\@vlinefalse%
                                        \else%
                                            \ifx#13%
                                                 \@nonefalse\@parenfalse\@brackfalse\@bracetrue\@vlinefalse%
                                            \else%
                                                 \ifx#14%
                                                       \@nonefalse\@parenfalse\@brackfalse\@bracefalse\@vlinetrue
                                                 \else%
                                                       \ifx#15%
                                                             \@nonefalse\@parenfalse\@brackfalse\@bracefalse\@vlinefalse%
                                                       \else%
                                                             \@nonetrue\@parenfalse\@brackfalse\@bracefalse\@vlinefalse%
                                                       \fi%
                                                 \fi%
                                            \fi%
                                        \fi%
                                   \fi%
                                   \if@smallmat%
                                        \if@none%
                                             \begin{smallmatrix}%
                                        \else%
                                            \if@paren%
                                                  \bigl(\begin{smallmatrix}%
                                            \else%
                                                  \if@brack%
                                                          \bigl[\begin{smallmatrix}%
                                                  \else%
                                                          \if@brace%
                                                               \bigl\{\begin{smallmatrix}%
                                                          \else%
                                                               \if@vline%
                                                                    \bigl\lvert\begin{smallmatrix}%
                                                                \else%
                                                                    \bigl\lVert\begin{smallmatrix}%
                                                                \fi%
                                                          \fi%
                                                  \fi%
                                            \fi%
                                        \fi%
                                   \else%
                                        \if@none%
                                             \begin{matrix}%
                                        \else%
                                            \if@paren%
                                                  \begin{pmatrix}%
                                            \else%
                                                  \if@brack%
                                                          \begin{bmatrix}%
                                                  \else%
                                                          \if@brace%
                                                               \begin{Bmatrix}%
                                                          \else%
                                                               \if@vline%
                                                                    \begin{vmatrix}%
                                                                \else%
                                                                    \begin{Vmatrix}%
                                                                \fi%
                                                          \fi%
                                                  \fi%
                                            \fi%
                                        \fi%
                                   \fi}%
                                  {\if@smallmat%
                                        \if@none%
                                             \end{smallmatrix}%
                                        \else%
                                            \if@paren%
                                                  \end{smallmatrix}\bigr)%
                                            \else%
                                                  \if@brack%
                                                          \end{smallmatrix}\bigr]%
                                                  \else%
                                                          \if@brace%
                                                               \end{smallmatrix}\bigr\}%
                                                          \else%
                                                               \if@vline%
                                                                    \end{smallmatrix}\bigr\rvert%
                                                                \else%
                                                                    \end{smallmatrix}\bigr\rVert%
                                                                \fi%
                                                          \fi%
                                                  \fi%
                                            \fi%
                                         \fi%
                                   \else%
                                        \if@none%
                                             \end{matrix}%
                                        \else%
                                            \if@paren%
                                                  \end{pmatrix}%
                                            \else%
                                                  \if@brack%
                                                          \end{bmatrix}%
                                                  \else%
                                                          \if@brace%
                                                               \end{Bmatrix}%
                                                          \else%
                                                               \if@vline%
                                                                    \end{vmatrix}%
                                                                \else%
                                                                    \end{Vmatrix}%
                                                                \fi%
                                                          \fi%
                                                  \fi%
                                            \fi%
                                        \fi%
                                   \fi}%
                         
\makeatother

\def\ger{\mathfrak}

\newcommand\CategoryTypeface{\mathbf}
\def\cat{\CategoryTypeface}

\newcommand\SheafTypeface{\mathcal}
\def\sh{\SheafTypeface}

\DeclareMathOperator\SPol{\cat{SPol}}

\DeclareMathOperator\Fun{\cat{Fun}}
\DeclareMathOperator\SDom{\cat{SDom}}
\DeclareMathOperator\SSp{\cat{SSp}}
\DeclareMathOperator\SMan{\cat{SMan}}

\DeclareScalableOperator{Sh}{\cat{Sh}}

\newcommand{\lBr}{[\kern-.65ex[}
\newcommand{\rBr}{]\kern-.65ex]}

\DeclareMathOperator\supp{\mathrm{supp}}

\DeclareMathOperator\Ob{\mathrm{Ob}}
\DeclareMathOperator\Spec{\mathrm{Spec}}

\DeclareScalableOperator{Ber}{\mathrm{Ber}}
\DeclareScalableOperator{Vol}{\mathrm{Vol}}
\DeclareScalableOperator{Form}{\Omega}
\DeclareScalableOperator{Ind}{\mathrm{Ind}}
\DeclareScalableOperator{Coind}{\mathrm{Coind}}
\DeclareScalableOperator{Der}{\mathrm{Der}}
\DeclareScalableOperator{Maps}{\mathrm{Maps}}
\DeclareScalableOperator{Ct}{\mathcal{C}} 
\DeclareScalableOperator{GCt}{\underline{\mathcal{C}}} 
\DeclareScalableOperator{Lp}{\mathbf{L}} 
\DeclareScalableOperator{Hom}{\mathrm{Hom}} 
\DeclareScalableOperator{GHom}{\underline{\mathrm{Hom}}} 
\DeclareScalableOperator{GEnd}{\underline{\mathrm{End}}} 
\DeclareScalableOperator{GPol}{\underline{\sh{P}}} 
\DeclareScalableOperator{GLin}{\underline{\sh{L}}} 
\DeclareScalableOperator{End}{\mathrm{End}} 
\DeclareScalableOperator{Aut}{\mathrm{Aut}} 
\DeclareScalableOperator{Uenv}{\mathfrak{U}} 
\DeclareScalableOperator{Cuenv}{\widehat{\mathfrak{U}}} 
\DeclareScalableOperator{ABer}{\Abs0{\mathrm{Ber}}}
\DeclareScalableOperator{Dist}{\mathcal{D}'}

\newcommand\Define[1]{\emph{#1}}

\begin{document}

\title{A convenient category of supermanifolds}
\author[A.~Alldridge]{Alexander Alldridge}

\address{Mathematisches Institut\\ Universit\"at zu K\"oln\\ Weyertal 86-90\\ 50931 K\"oln, Germany}
\email{alldridg@math.uni-koeln.de}
\thanks{This research was funded by the Leibniz independent junior research group grant and the SFB/Transregio 12 grant, provided by Deutsche Forschungsgemeinschaft (DFG), and the DFG--JSPS ``Lie Groups: Analysis and Geometry'' grant, provided by DFG and the Japan Society for the Promotion of Science (JSPS)}


\date{}

\begin{abstract}
	With a view towards applications in the theory of infinite-dimensio\-nal representations of finite-dimensional Lie supergroups, we introduce a new category of supermanifolds. In this category, supermanifolds of `maps' and `fields' (fibre bundle sections) exist. In particular, loop supergroups can be realised globally in this framework. It also provides a convenient setting for induced representations of supergroups, allowing for a version of Frobenius reciprocity. Finally, convolution algebras of finite-dimensional Lie supergroups are introduced and applied to a prove a supergroup Dixmier--Malliavin Theorem: The space of smooth vectors of a continuous representation of a supergroup pair equals the G\r{a}rding space given by the convolution with compactly supported smooth supergroup densities.
\end{abstract}


\maketitle

\section{Introduction}

Supermanifolds were defined in the 1970s by Berezin, Konstant, Leites, and Rogers, in order to provide a mathematically rigorous foundation for the supersymmetric field theories with bosonic and fermionic degrees of freedom. 

The well-established `Berezin--Kostant--Leites' approach is modeled on analytic and algebraic geometry, where one considers algebraic varieties or complex spaces as ringed spaces. So, one associates to an open set a ring of `functions' with `scalar values'. This ringed space approach allows nilpotents in the structural sheaves, so it is also convenient for the theory of (finite-dimensional) supermanifolds. 

However, the space of greatest interest to physics, notably, spaces of `maps'  or `fields' (\ie sections of fibre bundles), are infinite-dimensional. There have been several attempts to define an appropriate category of supermanifolds which contains such spaces, notably by Schmitt \cite{schmitt1} and Molotkov--Sachse \cite{molotkov,sachse,al-inf}. Remarkably, in the latter approach, Sachse--Wockel \cite{sachse-wockel} have proved the existence of a supergroup of superdiffeomorphisms for any compact supermanifold. 

However, the latter approach, while it is certainly the most general and versatible, is rather hard to handle, and in the former, the existence problem for supermanifolds of maps is not solved. Also, for  applications, it would be desirable to have a more concrete model at hand, which is closer to the ringed space approach in that supermanifolds should be a \emph{full} subcategory of a certain well-behaved larger category which possibly also contains not too distant relatives of smooth supermanifolds, such a (at least certain) superschemes, or analytic supermanifolds. As it turns out, this is easier to accomplish than one may think.

In fact, the main reason why the ringed space approach does not work well in infinite dimensions is that it is insufficient to consider \emph{scalar} valued functions on an infinite-dimensional space. Rather, one has to take locally defined maps with values in (somewhat) \emph{arbitrary} infinite-dimensional spaces into account. (In fact, values in all local models for the spaces one has in mind should suffice.) This was already observed by Douady in his thesis \cite{douady-probmodules}; he proposed the framework of \emph{espaces fonct\'es}\footnote{Neither the terminology nor the concept seem to have caught on, despite their naturality.} (functored spaces) which we will use below, and applied it in the solution of the moduli problem for compact analytic subspaces. 

As we show, this simple idea leads to a category of superspaces and and full subcategory of supermanifolds, in which the usual category of finite-dimensional supermanifolds is embedded fully faithfully. Moreover, supermanifolds of maps (\viz inner homs) and supermanifolds of fields (\viz sections of fibre bundles) exist in this category (\thmref{Th}{innerhom-smf} and \thmref{Prop}{inner-sections}). The formalism of Weil functors also works well in this setting, which allows a systematic treatment of tangent objects. In particular, the Lie superalgebra of a (possibly infinite-dimensional) Lie supergroup can be defined. 

We illustrate the utility of the new category of supermanifolds by three classes of examples. First, we give a brief account of loop and superloop supergroups, which can be treated with ease which the aid of inner hom functors. However, our main motivation comes from the theory of inifinite-dimensional representations of finite-dimensional supergroups. Therefore, secondly, we give an account of the basics of induced representations, proving Frobenius reciprocity (\thmref{Prop}{frobenius}). Thirdly, we show that natural convolution algebras of compactly supported distributions and compactly supported smooth densities can be defined for any finite-dimensional Lie supergroup $G$. Their action in $G$-representations (of arbitrary dimension) can be used to study questions of multiplicity (\thmref{Prop}{multi}). Moreover, we prove a super version of the Theorem of Dixmier--Malliavon which characterises the smooth vectors in a continuous supergroup pair representation as the G\r{a}rding space given by convolution with compactly supported smooth densities (\thmref{Th}{super-dixmal}). 

\section{An extension of the category of supermanifolds}

\subsection{Functored spaces}

\begin{Def}
	Let $\cat C$ be a category. A \Define{$\cat C$-functored space} is a pair $X=(X_0,\sh O_X)$ where $X_0$ is a topological space and $\sh O_X$ is a functor $\cat C\to\Sh0X$. We will write
	\[
		\sh O_X(U,A)=[\sh O_X(A)](U)\mathfa A\in\Ob\cat C\,,\,U\subset X\text{ open,}
	\]
	and $X^A=(X_0,\sh O_X(-,A))$. The functor $\sh O_X$ is called \Define{structure functor}.
	
	A \Define{morphism of $\cat C$-functored spaces} $f:X\to Y$ is $(f_0,f^*)$ where $f_0:X_0\to Y_0$ is a continuous map and $f^*:f_0^{-1}\sh O_Y\to\sh O_X$  a natural transformation. Here, $f_0^{-1}:\Sh0Y\to\Sh0X$ is the inverse image functor. With the obvious composition and identity morphisms, one obtains thus the category $\Fun_\cat C$ of $\cat C$-functored spaces. 
	
	If $X,Y$ are $\cat C$-functored spaces, then $Y$ is called an \Define{open subspace} of $X$ if, for all $A\in\Ob\cat C$, $Y^A$ is an open subspace of $X^A$. That is, $Y_0\subset X_0$ is open, and $\sh O_Y(-,A)=\sh O_X(-,A)|_{Y_0}$ \fa $A$. We write $Y=X_{Y_0}$. 
\end{Def}

We shall now construct categories $\cat C$ of `local models' for infinite-dimensional superspaces. This will allow us to apply the framework of functored spaces to the theory of superspaces, as we shall presently see. 

\begin{Par}
	A \Define{locally convex super-vector space} is by definition a $\ints/2\ints$ graded vector space $E=E_0\oplus E_1$ endowed with a locally convex Hausdorff vector space topology such that the graded parts are closed. Let $E,F$ be locally convex super-vector spaces and $U\subset E_0$ an open set. Let $\Ct[^\infty]0{U;F}$ denote the set of smooth maps $U\to F$. This is naturally a super-vector space. 
	
	For any $e\in E_0$ and $f\in\Ct[^\infty]0{U,F}$, let $\partial_ef$ be the derivative of $f$ in direction $e$. This defines an action of the symmetric algebra $S(E_0)$ on $\Ct[^\infty]0{U;F}$. Let 
	\[
		\Hom[_{S(E_0)}]0{S(E),\Ct[^\infty]0{U,F}}=\Set1{f:S(E)\to\Ct[^\infty]0{U,F}}{f\text{ even and $S(E_0)$-linear}}\ ,
	\]
	where $S(E)$ is the supersymmetric algebra.
	
	We write $f(P;x)=f(P)(x)$ for $f\in\Hom[_{S(E_0)}]0{S(E),\Ct[^\infty]0{U;F}}$, $P\in S(E)$, and $x\in U$. Then the $S(E_0)$-linearity of $f$ is expressed by 
	\[
		f(eP;x)=\partial_e[f(P)](x)\mathfa P\in S(E)\,,\,e\in E_0\,,\,x\in U\ ,
	\]
	the symbol $\partial_e$ denoting partial derivative in the direction of $e$.
	
	We let $\Ct[^\infty]0{E_U,F}$ be the set of all $f\in\Hom[_{S(E_0)}]0{S(E),\Ct[^\infty]0{U,F}}$ \scth for any $n\in\nats$, the map
	\[
		U\times E^n\to F:(x,e_1,\dotsc,e_n)\mapsto f(e_1\dotsm e_n;x)
	\]
	is smooth. Given any such $f$, we define $f_0=f(1)\in\Ct[^\infty]0{U,F}$. For any open set $V\subset F_0$, we let $\Ct[^\infty]0{E_U,F_V}$ be the subset of all $f\in\Ct[^\infty]0{E_U,F}$ \scth $f_0(U)\subset V$. 
	
	Assume given locally convex super-vector spaces $E,F,G$, open subsets $U\subset E_0$, $V\subset F_0$, $W\subset G_0$, and $f\in\Ct[^\infty]0{E_U,F_V}$, $g\in\Ct[^\infty]0{F_V,G_W}$. We wish to define the composite $g\circ f\in\Ct[^\infty]0{E_U,G_W}$. 
	
	To that end, consider the coproduct $\Delta:S(E)\to S(E)\otimes S(E)$, defined to be the unique even unital algebra morphism \scth $\Delta(e)=e\otimes1+1\otimes e$. Inductively, define $\Delta^k:S(E)\to S(E)^{\otimes(k+1)}$ by $\Delta^{k+1}=(\Delta\otimes\id^{\otimes k})\circ\Delta$. Similarly, let $\mu^k:S(E)^{\otimes(k+1)}\to S(E)$ be the $k$-fold product. We also set $\Delta^0=\mu^0=\id_{S(E)}$. 
	
	For homogeneous $e_1,\dotsc,e_n\in E$, we have
	\[
		\Delta^k(e_1\dotsm e_n)=\sum_{\mathbf n=\coprod_{j=0}^kI_j}\eps_{I_0,\dotsc,I_k}^{e_1,\dotsc,e_n}\cdot e_{I_0}\otimes\dotsm\otimes e_{I_k}\ .
	\]	
	Here, $\mathbf n=\{1,\dotsc,n\}$, $e_I=e_{i_1}\dotsm e_{i_m}$ for $I=\{i_1<\dotsm<i_m\}$ and $\eps_{I_0,\dotsc,I_k}^{e_1,\dotsc,e_n}=(-1)^N$,
	\[
		N=\#\Set1{i<j}{i\in I_p\,,\,j\in I_q\,,\,p>q\,,\,\Abs0{e_i}\equiv\Abs0{e_j}\equiv1\ (2)}\ .
	\]
	
	In the above formula for $\Delta^k$, the sum was over partitions $I_0,\dotsc,I_k$ of $\mathbf n$ where $I_j$ was possibly vacuous. Let $\Delta^{\prime k}$ be defined by restricting the summation to  those partitions where all $I_j\neq\vvoid$. 
	
	We define $f_{(n)}:S(E)\to\Maps0{U,S(F)}$ by $f_{(0)}=f$, and by 
	\[
		(n+1)!\cdot f_{(n)}(-;x)=\mu^n\circ(f(-;x))^{\otimes(n+1)}\circ\Delta^{\prime n}\mathfa n\sge1\,,\,x\in U\ .
	\]
	Thus, we define $g\circ f:S(E)\to\Ct[^\infty]0{E_U;G}$ by setting $(g\circ f)(1;x)=g_0(f_0(x))$ and 
	\begin{equation}\label{eq:super-faadibruno}
		(g\circ f)(P;x)=\sum_{n=0}^\infty g\Parens1{f_{(n)}(P;x);f_0(x)}\mathfa P\in S^k(E)\,,\,k\sge1\,,\,x\in U\ .
	\end{equation}
	This is a graded multivariate Fa\`a di Bruno formula.  
\end{Par}

\begin{Prop}[faadibruno-cat]
	Equation \eqref{eq:super-faadibruno} defines $g\circ f\in\Ct[^\infty]0{E_U,G_W}$. With this composition and the identity $\id=\id_{E_U}$ defined by 
	\[
		\id(1;x)=x\ ,\ \id(v;x)=v\ ,\ \id(P;x)=0
	\]
	\fa $x\in U$, $v\in E$, $P\in S^k(E)$, $k>1$, the pairs $(U,E)$, $(V,F)$ and the morphism sets $\Ct[^\infty]0{E_U,F_V}$ form a category. We denote $E_U=(U,E)$, \etc
\end{Prop}

\begin{proof}
	We have for $v\in E_0$, homogeneous $e_1,\dotsc,e_k\in E$, and $1\sle n<k$,
	\begin{align*}
		\Delta^{\prime n}(ve_1\dotsm e_k)=&\sum_{j=0}^n\sum_{\mathbf k=\coprod_{\ell=1}^nI_\ell\,,\,I_\ell\neq\vvoid}\eps^{e_1,\dotsc,e_k}_{I_1,\dotsc,I_n}\cdot e_{I_1}\otimes\dotsm\otimes e_{I_j}\otimes v\otimes e_{I_{j+1}}\otimes e_{I_n}\\
		&+\sum_{j=0}^n\sum_{\mathbf k=\coprod_{\ell=0}^nI_\ell\,,\,I_\ell\neq\vvoid}\eps^{e_1,\dotsc,e_k}_{I_0,\dotsc,I_n}\cdot e_{I_0}\otimes\dotsm\otimes ve_{I_j}\otimes\dotsm e_{I_n}\ .
	\end{align*}
	
	The definition of $f_{(n)}$ to the first summand gives $f(v;x)f_{(n-1)}(e_1\dotsm e_k;x)$; by applying it to the second, we obtain $[\partial_vf_{(n)}(e_1\dotsm e_k)](x)$, so 
	\[
		f_{(n)}(ve_1\dotsm e_k;x)=f(v;x)f_{(n-1)}(e_1\dotsm e_k;x)+[\partial_vf_{(n)}(e_1\dotsm e_k)](x)\ .
	\]
	 In fact, this statement remains true for $n=k$ (where the second term vanishes). 
	
	Thus, we obtain, for $P\in S^k(E)$, $k\sge1$,
	\begin{align*}
		[\partial_v(g&\circ f)(P)](x)=\sum_{n=0}^{k-1}[\partial_{f(v;x)}g(f_{(n)}(P;x))](f_0(x)+g\Parens1{[\partial_vf_{(n)}(P)](x);f_0(x)}\\
		&=\begin{aligned}[t]
			&g\Parens1{[\partial_vf_{(0)}(P)](x);f_0(x)}+g\Parens1{f(v;x)f_{(k-1)}(P;x);f_0(x)}\\
			&+\sum_{n=1}^{k-1}g\bigl(f(v;x)f_{(n-1)}(P;x)+[\partial_vf_{(n)}(P)](x);f_0(x)\bigr)
		     \end{aligned}\\
		    &=\sum_{n=0}^kg\Parens1{f_{(n)}(vP;x);f_0(x)}=(g\circ f)(vP;x)\ .
	\end{align*}
	This equation also hold trivially for $k=0$. 
	
	Therefore, $g\circ f\in\Hom[_{S(E_0)}]0{S(E);\Ct[^\infty]0{U,F}}$. In fact $g\circ f\in\Ct[^\infty]0{E_U,F_V}$, as is easy to check. The associativity of $\circ$ is an easy consequence of the coassociativity equation $(\Delta\otimes\id)\circ\Delta=(\id\otimes\Delta)\circ\Delta$, which implies the same statement for $\Delta'$. 
	
	What remains to be shown is that $\id$ acts as the identity for the composition. To that end, observe that for $1\sle n\sle k$, 
	\[
		\id_{(n)}(e_1\dotsm e_k;x)=\delta_{n+1,k}\cdot e_1\dotsm e_k\ .
	\]
	This readily implies the assertion. 
\end{proof}

\begin{Rem}
	If $E,F$ are f.d.~super-vector spaces and $U\subset E$ is open, then 
	\[
		\Ct[^\infty]0{E_U,F}=\Parens1{\Hom[_{S(E_0)}]0{S(E),\Ct[^\infty]0{U,\reals}}\otimes F}_0=\Parens1{\Ct[^\infty]0{U,\reals}\otimes\textstyle\bigwedge E_1^*\otimes F}_0\ ,
	\]
	since $S(E)=S(E_0)\otimes\bigwedge E_1$. This motivates the following definition. 
\end{Rem}

\begin{Def}
	We denote the category defined in \thmref{Prop}{faadibruno-cat} by $\SDom$ and call this the category of \Define{smooth superdomains}. 
\end{Def}

\begin{Rem}
	The same construction can be performed for real-analytic and holomorphic maps (working over $\cplxs$ in the latter case); \vq\cite{bochnak-siciak,BGlN} for the correct definition of analyticity in this setup.
	
	Working over $\reals$, one may take $\Ct[^\varpi]0{U,F_\cplxs}$ (where $\varpi=\infty,\omega$) and complex-linear $\Hom[_{S(E_{0\cplxs})}]0{S(E_\cplxs),\cdot}$ in the definition of hom-sets, to obtain categories which correspond to possibly infinite-dimensional versions of what Bernstein \cite{deligne-morgan} calls (smooth resp.~analytic) \emph{cs} manifolds. 
\end{Rem}

\begin{Def}
	Let $E$ be a locally convex super-vector space. Define a $\SDom$-functored space $L(E)=(E_0,\sh O_E)$ where
	\[
		\mathcal O_E(U,F_V)=\Ct[\infty]0{E_U,F_V}
	\]
	and for $g:F_V\to(V',F )$, 
	\[
		\sh O_E(U,g):\sh O_E(U,F_V)\to\sh O_E(U,F'_{V'}):f\mapsto g\circ f\ .
	\]
	This is a sheaf with restriction maps $f\mapsto f|_{U'}$ defined by 
	\[
		f|_{U'}(P;x)=f(P;x)\mathfa P\in S(E)\,,\,x\in U'\subset U\ .
	\]
	
	We call the functored space $L(E)$ \Define{linear}. Open subspaces of such linear functored spaces will be called \Define{smooth superdomains}. 
\end{Def}

\begin{Prop}[mor-fun]
	Let $X$ be a $\SDom$-functored space and $F_V\in\Ob\SDom$. Assume that $X$ has an open cover by smooth superdomains. 
	
	Then we have a natural isomorphism
	\[
		\phi:\Hom[_{\Fun_{\SDom}}]0{X,L(F)_V}\to\sh O_X(X_0,F_V)\ .
	\]
	In particular, $E_U\mapsto L(E)_U$ is a fully faithful embedding of $\SDom$ in $\Fun_{\SDom}$. 
\end{Prop}

Douady \cite{douady-probmodules} comments on the proof of a version of the proposition: ``La d\' emonstra\-tion est asinitrottante.'' We feel inclined to endorse this view. Nevertheless, we give the relevant construction. 

\begin{proof}
	The proposition is in fact quite more general and holds for any category $\cat C$ which is enriched over sheaves on topological spaces. We do not give the most general formulation here. The proof is similar to that of the Yoneda lemma.\footnote{In fact, it is an instance of the enriched Yoneda lemma \cite{kelly-basicenrichedcat}.}
	
	Indeed, let us define the map for $X=L(E)_U$. This will be sufficient, since by naturality we will then be able to glue it to a globally defined map. In the special case, let $f:L(E)_U\to L(F)_V$ be a morphism of functored spaces. Then
	\[
		\Gamma(f^*_{F_V}):\Ct[^\infty]0{F_V,F_V}\to\Ct[^\infty]0{E_U,F_V}\ ,
	\]
	so we set $\phi(f)=\Gamma(f^*_{F_V})(\id_{F_V})\in\Ct[^\infty]0{E_U,F_V}=\sh O_{L(E)_U}(U,F_V)$. 
	
	We may define a converse map $\psi$ by setting for $f\in\Ct[^\infty_{E,F}]0{U,V}$, $\psi(f)_0=f_0$, and defining for $G_W\in\Ob\SDom$ a sheaf morphism 
	\[
		\psi(f)^*_{G_W}:\sh O_{L(F)_V}(-;G_W)\to f_{0*}\sh O_{L(E)_U}(-;G_W)
	\]
	by setting, for any open $V'\subset V$ and $h\in\Ct[^\infty]0{F_{V'},G_W}$, 
	\[
		\psi(f)^*_{G_W,V'}=h\circ f|_{f_0^{-1}(V')}\in\Ct[^\infty]0{F_{f_0^{-1}(V')},G_W}=f_{0*}\sh O_{L(E)_U}(-;G_W)(V')\ .
	\]
	
	It is obvious that $\phi\circ\psi=\id$, and $\psi\circ\phi=\id$ follows by using the naturality of morphisms of functored spaces. Hence follows our claim. 
\end{proof}

\begin{Par}
	For locally convex super-vector spaces $E,F$, let $\sh P(E,F)\subset\Hom0{S(E),F}$ consist of the elements $f$ \scth
	\[
		(e_1,\dotsc,e_k)\mapsto f(e_1\dotsm e_k):E^k\to F
	\]
	is continuous for any $k\in\nats$. The elements of $\sh P(E,F)$ are called \Define{polynomial maps}. 
	
	If $E=E_0$, $F=F_0$, then this notion coincides with the classical definition of polynomials maps, \cf\cite{bochnak-siciak-polynomial}. To any $f\in\sh P(E,F)$, we may associate an element $\tilde f$ of $\Hom[_{S(E_0)}]0{S(E),\sh P(E_0,F_0)\oplus\sh P(E_0,\Pi F_1)}$ by 
	\[
		\tilde f(P;x)=\sum_{k=0}^\infty\frac1{k!} f(x^kP)\mathfa P\in S(E)\,,\,x\in E .
	\]
	(Here, $\Pi$ denotes the reversal of parity functor.)
	
	This defines an injection $\sh P(E,F)\to\Ct[^\infty_E]0{E_0,F}$, by which we identify $\sh P(E,F)$ with its image.
	If $U\subset E_0$ is open, let $\sh P_E(U,F)$ be the set of the restrictions $f|_U$, for $f\in\sh P(E,F)$. For any open $V\subset F_0$, let $\sh P_{E,F}(U,V)=\sh P_E(U,F)\cap\Ct[^\infty_{E,F}]0{U,V}$. The following lemma is immediate.
\end{Par}

\begin{Lem}
	Taking as objects pairs $E_U$, $F_V$ where $E,F$ are locally convex super-vector spaces and $U\subset E$, $V\subset F$ are open, and $\sh P_{E,F}(U,V)$ as morphism sets, the structure induced by $\SDom_{\sh C^\infty}$ defines a subcategory $\SPol$. 
\end{Lem}

\begin{Lem}[pol-emb]
	The functor $\Fun_{\SDom}\to\Fun_{\SPol}$ given by the restriction of functored spaces to the subcategory of algebraic superdomains is fully faithful. 
\end{Lem}

\begin{proof}
	This is trivial, since the index categories the functored spaces are defined on have the same classes of objects.
\end{proof}

\begin{Def}
	Let $\cat C=\SPol$. An $X\in\Ob\Fun_{\cat C}$ is called \Define{unital} if endowed with embeddings $V\subset\sh O_X(U,F_V)$ \fa open $\vvoid\neq U\subset X_0$, $F_V\in\cat C$, compatible with the presheaf structure. A morphism of unital functored spaces is called \Define{unital} if it respects these embeddings.  
	
	Let $\SSp$ denote the category whose objects are the unital functored spaces and whose morphisms are the unital morphisms thereof. The objects of $\SSp$ are called \Define{superspaces}. 
	
	The superdomain $L(E)_U\in\Fun_{\SDom}$ induces a superspace denoted $E_U\in\SSp$. 
	The unitality condition for morphisms of functored spaces between superdomains is automatically verified, so the following statement is immediate from \thmref{Prop}{mor-fun} and \thmref{Lem}{pol-emb}. 
\end{Def}

\begin{Cor}[super-mor]
	Let $X\in\SSp$ be covered by smooth superdomains and $F_V$ a smooth superdomain. Then there is a natural bijection
	\[
		\Hom[_{\SSp}]0{X,F_V}\cong\sh O_X(X_0,F_V)\ .
	\]
\end{Cor}

\begin{Rem}
	From what has been proved above, the categories of f.d.~smooth supermanifolds have a fully faithful embedding into $\SSp$. Similarly, f.d.~real analytic supermanifolds embed into this category. Moreover, f.d.~smooth and analytic \emph{cs}Êmanifolds and complex analytic supermanifolds embed fully faithfully into a similar category of complex superspaces. 
	
	We have the following pasting lemma.
\end{Rem}

\begin{Lem}[super-glue]
	Assume given superspaces $X_i$, open subspaces $X_{ij}\subset X_i$, and isomorphisms $\psi_{ij}:X_{ji}\to X_{ij}$ \scth $\psi_{ij}\circ\psi_{jk}=\psi_{ik}$. Then there exists a superspace $X$ with open subspaces $U_i\subset X$ and isomorphisms $\psi_i:X_i\to U_i$ \scth $\psi_i\circ\psi_{ij}=\psi_j$. It is determined uniquely up to unique ismorphism by the following universal property: For any superspace $Y$ and any morphisms $\phi_i:X_i\to Y$ \scth $\phi_i\circ\psi_{ij}=\phi_j$, there exists a unique morphism $\phi:X\to Y$ \scth $\phi\circ\psi_i=\phi_i$. 
\end{Lem}

\begin{Def}
	A \Define{non-Hausdorff supermanifold} is a superspace $X$ with an open cover by superdomains. If the open cover may be chosen from open subspaces of some class $\sh C$ of locally convex super-vector spaces, then we say that $X$ is locally modelled over $E\in\sh C$; in particular, this defines the notion of $X$ being locally of finite dimension. A \Define{supermanifold} is a non-Hausdorff supermanifold $X$ whose underlying topological space $X_0$ is Hausdorff. 
	
	We thus obtain a full subcategories $\SMan$ and $\SMan_{NH}$ of $\SSp$ consisting of supermanifolds resp.~non-Hausdorff supermanifolds. The pasting lemma carries over to non-Hausdorff supermanifolds, and to supermanifolds under conditions which insure the Hausdorff property for the topological space underlying the glued superspace.
\end{Def}

\begin{Cor}[super-prod]
	The category $\SSp$ has a terminal object. Finite products exist in the subcategories of non-Hausdorff supermanifolds and supermanifolds. 
\end{Cor}

\begin{proof}
	The latter statement in the case of binary products is local up to the Hausdorff condition, and thus follows in the usual way from \thmref{Cor}{super-mor}. As to the former, let $*$ be the functored space on the point with structure functor defined by $\sh O_*(*,F_V)=V$. The claim is now an obvious restatement of the unitality of morphisms of superspaces. 
\end{proof}

\begin{Rem}
	The embedding of finite-dimensional supermanifolds into $\SMan$ preserves finite products. Indeed, it is sufficient to check this for superdomains, where is true by construction. 
\end{Rem}

\subsection{Supermanifolds of maps}

We will now show that under very mild conditions, inner homs (`supermanifolds of maps') exist in the category of supermanifolds. We first recall some known facts concerning the existence of inner homs in the category of Hausdorff topological spaces.

\begin{Par}
	Let $X$ be a Hausdorff space. A subset $A\subset X$ is \Define{$k$-closed} (resp.~$k$-open) if $K\cap A$ is closed (resp.~open) in $K$ for any compact $K\subset X$. The set of $k$-open subsets of $X$ forms a topology, the \Define{$k$-topology}, which is finer than $X$. Let $kX$ be $X$ with this topology. Then $X$ is a \emph{$k$-space} if $X=kX$ as topological spaces. Observe that $kX$ is a $k$-space. Open and closed subspaces of $k$-spaces are $k$-spaces. If $X$ is sequential (\ie every sequentially closed subset is closed) then $X$ is a $k$-space. In particular, this holds for first-countable spaces. If $Y$ is Hausdorff space \scth $X$ and $Y$ are first-countable or $X$ is a $k$-space and $Y$ is locally compact, then $X\times Y$ is a $k$-space. If $X$ is a $k$-space, then $f:X\to Y$ is continuous if and only if $f:X\to kY$ is.  
	
	If $X$ and $Y$ are Hausdorff spaces, and endow $\Ct0{X,Y}$ with the compact-open topology; a subbase for the topology is given by the sets $W_{K,U}$ of maps $f:X\to Y$ \scth $f(K)\subset U$, for all compact $K\subset X$ and open $U\subset Y$. If $Y$ is a locally convex vector space, then $\Ct0{X,Y}$ is a locally convex vector space. If in addition, $X$ is locally compact and metrisable, and $Y$ is metrisable, then $\Ct0{X,Y}$ is metrisable. The following lemma is well-known.
\end{Par}

\begin{Lem}[curry]
	Let $X,Y,Z$ be Hausdorff spaces \scth $X\times Y$ and either $X$ or $\Ct0{Y,Z}$ are $k$-spaces. Then 
	\[
		\Ct0{X,\Ct0{Y,Z}}\to\Ct0{X\times Y,Z}:f\mapsto((x,y)\mapsto f(x)(y))
	\]
	is a bijection. In fact, if $X\times Y$ is replaced by $k(X\times Y)$, and $\Ct0{Y,Z}$ by $k\Ct0{Y,Z}$, then this holds without any restriction on $X,Y,Z$ (besides being Hausdorff).
\end{Lem}

\begin{Par}[smooth-top]
	Let $E$, $F$ be locally convex super-vector spaces and $U\subset E_0$ be open. We have an injection
	\[
		\Ct[^\infty_E]0{U,F}\to\prod\nolimits_{k,\ell}\Ct0{U\times E_0^k\times E_1^\ell,F}:f\mapsto((x,e_1,\dotsc, e_{k+\ell})\mapsto f(e_1\dotsm e_{k+\ell};x))\ .
	\]
	
	We endow $\Ct[^\infty_E]0{U,F}$ with the subspace topology induced by this injection. Therefore, if $F$ is locally convex and metrisable, then $\Ct[^\infty_E]0{U,F}$ is locally convex. If in addition, $E$ is finite dimensional and $F$ is metrisable, then $\Ct[^\infty_E]0{U,F}$ is metrisable. 
\end{Par}

\begin{Def}
	We will consider certain subclasses of locally convex super-vector spaces, namely, the classes $\cat{fd}$, $\cat{fc}$, $\cat{met}$, of, respectively, finite-dimensional, first countable and metrisable locally convex super-vector spaces. 
\end{Def}

\begin{Lem}[smooth-curry]
	Let $E,F,G$ be locally convex vector spaces where $(E,F)\in\cat{fc}\times\cat{fc}$. For any open subsets $U\subset E$, $V\subset F$, and $W\subset G$, there is a natural bijection
	\[
		\Ct[^\infty]0{U,\Ct[^\infty]0{V,W}}\to\Ct[^\infty]0{U\times V,W}
	\]
\end{Lem}

\begin{proof}
	Consider the evaluation map $e:\Ct[^\infty]0{V,W}\times V\to W:(f,v)\mapsto f(v)$. Inductively, one shows that it possesses an $n$-th derivative
	\[	
		d^ne(f,x)(h_1,v_1,\dotsc,h_n,v_n)=d^nf(x)(v_1,\dotsc,v_n)+\sum\nolimits_id^{n-1}h_i(v_1,\dotsc,\widehat{v_i},\dotsc,v_n)
	\]
	which is continuous in view of \thmref{Lem}{curry}. (The lemma applies since $U$ and $U\times V$ are $k$-spaces if $E$ and $F$ are first-countable.) Since $W$ is locally convex, $e$ is smooth (\vq\cite[Proposition 7.4]{BGlN}). Thus, if $f\in\Ct[^\infty]0{U,\Ct[^\infty]0{V,W}}$, then the map $g:U\times V\to W$, defined by $g(x,y)=f(x)(y)$, is smooth, since $g=e\circ(f\times\id)$. 
	
	Conversely, let $g\in\Ct[^\infty]0{U\times V,W}$, and define $f(x)(y)=g(x,y)$. Then we have $f:U\to\Ct[^\infty]0{V,W}$. Moreover, for $x\in U$, $y\in V$, $d^n(f(x))(y)=d^n_2g(x,y)$ where $d_j$ denotes derivatives with respect to the $j$th argument. Since $d_2g$ is continuous, \thmref{Lem}{curry} implies that $x\mapsto d^n(f(x)):U\to\Ct0{V\times F^n,G}$ is continuous, so that $f:U\to\Ct[^\infty]0{V,W}$ is continuous. Inductively, its $k$-th derivative is given by $d^kf(x)(v_1,\dotsc,v_k)(y)=d_1^kg(x,y)(v_1,\dotsc,v_k)$, and this is continuous as a map $U\times E^k\to\Ct[^\infty]0{V,G}$ by a similar argument. Since the latter space is locally convex, the assertion is proved by applying \loccit\ again.
\end{proof}

\begin{Rem}
	An avenue to further generalisation seems feasible in the framework of linear $k$-spaces (\vq Ref.~\cite{froelicher-jarchow-kvects}). To that end, one is lead to investigate when $\Ct[^\infty]0{V,W}$ is a $ck$-space (in the terminology of \opcit). This might be interesting since bornological spaces (and in particular, LF spaces) fall in this category. 
\end{Rem}

\begin{Par}
	For super-vector spaces $E,F$, let $\GHom0{E,F}$ denote the set of \emph{all} linear maps. Then $\GHom0{E,F}=\Hom0{E,F\oplus\Pi F}$. In particular, for $E_U,F_V\in\SDom_{\sh C^\infty}$, we have
	\[
		\GHom[_{S(E_0)}]0{S(E),\Ct[^\infty]0{U,F}}=\Hom[_{S(E_0)}]0{S(E),\Ct[^\infty]0{U,F\oplus\Pi F}}\ .
	\]
	We therefore define
	\[
		\GCt[^\infty]0{E_U,F}=\Ct[^\infty]0{E_U,F\oplus\Pi F}\ .
	\]
\end{Par}

\begin{Lem}[super-curry]
	Let $E_U,F_V,G_W\in\Ob\SDom_{\sh C^\infty}$ where $(E,F)\in\cat{fc}\times\cat{fc}$. There is a bijection
	\[
		\Ct[^\infty]0{E_U,\GCt[^\infty]0{F_V,G}_{\Ct[_{F,G}^\infty]0{V,W}}}\to\Ct[^\infty]0{(E\times F)_{U\times V},W_G}
	\]
	which is natural in $E_U$, considered as an object of $\SSp$. 
\end{Lem}

\begin{proof}
	Observe that the statement of the lemma is meaningful, because $\Ct[^\infty]0{F_V,G_W}$ is an open subset of $\GCt[^\infty]0{F_V,G}_0$. Note that $S(E\times F)=S(E)\otimes S(F)$. We define  
	\[
		\phi:\Hom1{S(E),\Ct[^\infty]0{U,\Ct[^\infty]0{F_V,G_W}}}\to\Hom1{S(E)\otimes S(F),\Ct[^\infty]0{U\times V,G}}
	\]
	by the assignment 
	\[
		\phi(f)(P\otimes Q;x,y)=[f(P;x)](Q;y)\ .
	\]
	
	The parity of the right hand side is $\Abs0P+\Abs0Q$ in the case of homogeneous $P,Q$, so that $\phi(f)$ is indeed even.  It is easy to check that $f$ is $S(E_0)$-linear if and only if $\phi(f)$ is $S(E_0)\otimes S(F_0)$-linear (the $\reals\otimes S(F_0)$-linearity is automatic, since it is verified by $f(P;x)$ for any $P,x$). It now follows from \thmref{Lem}{smooth-curry} that 
	\[
		(x,y,e_1,\dotsc,e_n,f_1,\dotsc,f_m)\mapsto \phi(f)(e_1\dotsm e_n\otimes f_1\dotsm f_m;x,y)
	\]
	is a smooth map $U\times V\times E^n\times F^m\to G$ \fa $m$ if and only if 
	\[
		U\times E^n\to\GCt[^\infty]0{F_V,G}:(x,e_1,\dotsc,e_n)\mapsto f(e_1\dotsm e_n;x)
	\]
	is smooth. The naturality is easy to check using \eqref{eq:super-faadibruno}. This implies the claim.
\end{proof}

We can now show the existence of inner homs. We recall some useful terminology. 

\begin{Def}
	Given any category $\cat C$ and $S,X\in\Ob\cat C$, we denote by $X(S)$ the set of morphisms $S\to X$. Instead of $x\in X(S)$ we write $x\in_SX$, and call $x$ an \Define{$S$-point of $X$}. For $f:X\to Y$ and $x\in_SX$, we define $f(x)=f\circ x\in_SY$. 
	
	Then the Yoneda lemma can be stated as follows: Given, for fixed $X,Y\in\Ob\cat C$ and arbitrary $S\in\Ob\cat C$, set maps $f_S:X(S)\to Y(S)$, there is $f:X\to Y$ \scth $f(x)=f_S(x)$ \fa $x\in_SX$ if and only if 
	\[
		f_T(x\circ g)=f_S(x)\circ g\mathfa g:T\to S\,,\,x\in_SX\ .
	\]
	In this case, $f$ is unique. 

	A morphism $p:X\to S$ of superspaces is also called a \Define{superspace $X$ over $S$}, denoted $X/S$ and $p=p_X$. A \Define{morphism $f:X/S\to Y/S$ of superspaces over $S$} is a morphism $f:X\to Y$ of superspaces \scth $p_Y\circ f=p_X$. 
	
	Hence the category $\SSp_S$ of superspaces over $S$. We denote the hom sets in this category by $\Hom[_S]0{-,-}$. Observe that $\Hom[_S]0{S\times X,S\times Y}=\Hom0{S\times X,Y}$, where $f\in\Hom[_S]0{S\times X,S\times Y}$ and $g\in\Hom0{S\times X,Y}$ are related by the equation $f(s,x)=(s,g(s,x))$ \fa $s\in_TS$, $x\in_TX$. 
	
	If $\cat C$ is some class of locally convex super-vector spaces and $X\in\Ob\SSp$, then we say that $X$ is locally in $\cat C$ if $X$ has an open cover by subspaces $Y$ isomorphic to superdomains $E_U$ with $E\in\cat C$. 
\end{Def}

\begin{Th}[innerhom-smf]
	Let $\SMan_{NH}$ denote the category of smooth non-Hausdorff supermanifolds which are locally in $\cat{fc}$. Let $X,Y\in\Ob\SMan_{SH}$ where $X$ is locally in $\cat{fd}$ and $Y$ is locally in $\cat{met}$.  The functor defined for $T\in\Ob\SMan_{SH}$ by
	\[
		\GHom0{X,Y}(T)=\Hom0{T\times X,Y}=\Hom[_T]0{T\times X,T\times Y}
	\]
	is representable by a non-Hausdorff supermanifold which is locally in $\cat{met}$. Moreover, $\GHom0{X,Y}_0$ is Hausdorff whenever $Y_0$ is. 
\end{Th}

\begin{Lem}[hausdorff-cond]
	Let $(U_\alpha,U_{\alpha\beta},\psi_{\alpha\beta})$ be gluing data in topological spaces, where $U_\alpha$ are Hausdorff. Let $U$ be the space obtained by gluing these data. Then $U$ is Hausdorff if and only if the map
	\[
		\textstyle\coprod_{\alpha\beta}U_{\alpha\beta}\to\coprod_\alpha U_\alpha\times\coprod_\beta U_\beta
	\]
	given by $\coprod_{\alpha\beta}(\id,\psi_{\beta\alpha})$ has closed range.
\end{Lem}

\begin{proof}
	The range of the map in the statement of the lemma is the graph of the equivalence relation $R$ defining the glued space $U=U'/R$ where $U'=\coprod_\alpha U_\alpha$. By \cite[Chapter I, ¤5.2, Example 2) and ¤8.3, Proposition 8]{bourbaki-top}, the relation $R$ is open (\ie the canonical map $U'\to U$ is open), and $U$ is Hausdorff if and only if the graph of $R$ is closed in $U'$. 
\end{proof}

\begin{proof}[\protect{Proof of \thmref{Th}{innerhom-smf}}]
	For non-Hausdorff supermanifolds, the statement is entirely local, by the local definition of products of non-Hausdorff supermanifolds. By \thmref{Lem}{super-curry}, \thmref{Cor}{super-mor}, $\GHom0{E_U,F_V}$ is represented by $G_W$ where we set $G=\GCt[^\infty]0{E_U,F}$ and $W=\Ct[^\infty]0{E_U,F_V}$. It remains to prove that $\GHom0{X,Y}_0$ is Hausdorff whenever $Y_0$ is. 
	
	To that end, we observe first that if $A,B,C$ are non-Hausdorff supermanifolds, then there is a natural morphism
	\[
		\circ:\GHom0{A,B}\times\GHom0{B,C}\to\GHom0{A,C}
	\]
	whenever the inner homs exist. Indeed, for any $T$, $\circ$ is given on $T$-points by the composition in $\SSp_T$, 
	\[
		\circ:\Hom[_T]0{T\times A,T\times B}\times\Hom[_T]0{T\times B,T\times C}\to\Hom[_T]0{T\times A,T\times C}\ .
	\]
	Since $\Hom0{A,B}=\GHom0{A,B}_0$ as sets, we have that with the induced topology, the ordinary composition
	\[
		\Hom0{A,B}\times\Hom0{B,C}\to\Hom0{A,C}
	\]
	is continuous (as the underlying map of a morphism of superspaces). 
	
	Let $(X_i)$ be an open cover of $X$ by (finite-dimensional) superdomains. The morphisms $X_i\to X$ induce a monomorphism 
	\[
		\GHom0{X,Y}\to\textstyle\coprod\nolimits_i\GHom0{X_i,Y}
	\]
	which gives a continuous injection on the level of topological spaces. Hence, it is sufficient to prove that $\GHom0{X,Y}_0$ is Hausdorff in the case of a superdomain. If $Y$ is a superdomain, we already know that this is the case. 
	
	Let $(Y_i)$ be an open cover of $Y$ by (metrisable locally convex) superdomains. Then $H=\GHom0{X,Y}$ is the colimit of the diagram
	\[
		\xymatrix{%
			\coprod\nolimits_{ij}\GHom0{X,Y_{ij}}\ar@<0.5ex>[r]\ar@<-0.5ex>[r]&\coprod\nolimits_j\GHom0{X,Y_j}=H'}
	\]
	and $Y$ is the colimit of the diagram $\xymatrix@1{\coprod\nolimits_{ij}Y_{ij}\ar@<0.5ex>[r]\ar@<-0.5ex>[r]&\,\coprod\nolimits_jY_j=Y'}$.  
	
	By \thmref{Lem}{hausdorff-cond}, the graph of the equivalence relation $R_Y$ on $Y_0'$ \scth $Y_0=Y'_0/R_Y$ is closed in $Y_0'\times Y_0'$, and in order to prove that $H_0$ is Hausdorff, it is sufficient to show that the graph of the relation $R_H$ \scth $H_0=H'_0/R_H$ is closed. Hence, let $f,g\in H_0'$ be $R_H$-inequivalent.
	
	Let $X=E_U$, $Y_j=(F_j)_{V_j}$. Then $H_0'=\coprod_j\Hom0{X,Y_j}=\coprod_j\Ct[^\infty]0{E_U,F_{j,V_j}}$; there are open $V_{ij}\subset V_i$ \scth $Y_{ij}=(F_i)_{V_{ij}}$ and isomorphisms $\psi_{ij}:Y_{ji}\to Y_{ij}$ defining the gluing data. We may assume $f\in\Hom0{X,Y_i}$, $g\in\Hom0{X,Y_j}$. Then 
	\[
		f\notin\Hom0{X,Y_{ij}}\,,\,g\notin\Hom0{X,Y_{ji}}\text{ or }f\neq \psi_{ij}\circ g\ .
	\]
	
	In a first step, assume that $f_0$ is inequivalent to $g_0$ for the equivalence relation on $\coprod_j\Ct[^\infty]0{U,V_j}$ induced by the gluing data $(\Ct[^\infty]0{U,V_i},\Ct[^\infty]0{U,V_{ij}},\psi_{ij,0}\circ(-))$. Then 
	\[
		f_0(U)\not\subset V_{ij}\,,\,g_0(U)\not\subset V_{ji}\text{ or }f_0\neq \psi_{ij,0}\circ g_0\ .
	\]
	Hence, there exists $x\in U$ \scth 
	\[
		f_0(x)\notin V_{ij}\,,\,g_0(x)\notin V_{ji}\text{ or }f_0(x)\neq \psi_{ij,0}(g_0(x))\ .
	\]
	Since $R_Y$ has closed graph, there are open neighbourhoods $U_i\subset V_i$, $U_j\subset V_j$ of $f(x)$ resp.~$g(x)$, \scth for all $(a,b)\in U_i\times U_j$,
	\[
		a\notin V_{ij}\,,\,b\notin V_{ji}\text{ or }a\neq\psi_{ij,0}(b)\ .
	\]
	
	The sets $W_{x,U_k}\subset\Ct0{U,V_k}$, for $k=i,j$, are open neighbourhoods of $f_0$ resp.~$g_0$ in the compact-open topology, and no $(f',g')\in Y_0'\times Y_0'$ \scth $(f'_0,g'_0)\in W_{x,U_i}\times W_{x,U_2}$ intersects the graph of $R_H$. By the definition of the topology on $\Ct[^\infty]0{E_U,F_{V_k}}$, it follows that $(f,g)$ is not in the closure of $R_H$. 
	
	Next, we assume that $f_0$ is equivalent to $g_0$, so $f_0(U)\subset V_{ij}$, $g_0(U)\subset V_{ji}$, and $f_0=\psi_{ij,0}\circ g_0$. Then we may assume $0\in V_i$, $0\in V_j$, and that $f_0=g_0=0$. Hence, we have $f\neq\psi_{ij}\circ g$. Since $\Ct[^\infty]0{E_U,F_{V_i}}$ is Hausdorff, there are disjoint open subsets $U_i,U_j'\subset\Ct[^\infty_{E,F_i}]0{U,V_i}$ \scth $f\in U_i$, $\psi_{ij}\circ g\in U_j'$. By the continuity of composition, there is an open neighbourhood $U_j\subset\Ct[^\infty]0{E_U,F_{j,V_j}}$ of $g$ \scth $\psi_{ij}\circ h\in U_j'$ \fa $h\in U_j$. Then $U_i\times U_j$ is an open neighbourhood of $(f,g)$ in $Y_0'\times Y_0'$ which does not intersect the graph of $R_H$. This proves the claim. 
\end{proof}

\subsection{Bundles and supermanifolds of sections}

In what follows, for definiteness, we work in the categories $\SSp$ of real superspaces and $\SMan$ of real smooth supermanifolds. Most statements, remain true for the case of \emph{cs} and complex superspaces, and in the analytic case.

\begin{Def}
	If $X/S\in\SSp_S$ and $U\subset S$ is an open subspace, define the open subspace $X_U=p_X^{-1}(U)$ of $X$ to correspond to the open subset $p_{X0}^{-1}(U_0)\subset X_0$. 
	
	A supermanifold $X/S$ where $S$ is also a supermanifold is called a \Define{family} if $S$ has a cover $(U_i)$ by open subspaces \scth \fa $U=U_i$, there is some supermanifold $F$ \scth $p_2:F\times U\to U$ \fs $F$ is isomorphic to $X_U/U$ in $\SSp_U$. 
	
	Any such isomorphism $U\times V/U\to X_U/U$ is called a \Define{local trivialisation} of $X$ over $U$; we denote the set of such local trivialisations by $\tau_X(U)$ (for arbitrary $F$) resp.~$\tau_X^F(U)$ (for fixed $F$). Then $\tau_X$ and $\tau_X^F(U)$ are sheaves of sets on $S$; for example, $\tau_X^F(U)$ is the subset of isomorphisms in $\Hom[_U]0{F\times U,X_U}$. 
	
	If $(U_i)$ and $F$ may be chosen \scth each $U_i$ is contained in a connected component of $S$ and $V$ is independent of $i$ for all $U_i$ contained in a fixed connected component, then we call $X$ a \Define{bundle}. All of the above may also be defined for non-Hausdorff supermanifolds. 
\end{Def}

\begin{Rem}
	There is an obvious pasting lemma for families and for bundles. In the case of bundles, the Hausdorff condition for the glued bundle is automatic if it is verified by the gluing data. In particular, if $X/T$ is a bundle and $f:S\to T$ is a morphism, then $S\times_TX$ exists and is a bundle over $S$. As a very special case, if $t\in T_0$, the canonical morphism $*\to T$ given by $t$ is a supermanifold over $T$. We set $X_t=t\times_TX$ and call this the \Define{fibre of $X$ over $t$}. 
\end{Rem}

\begin{Prop}[inner-sections]
	Let $X/S$ be a bundle with fibre $F$. The set-valued functor on $\SMan$ defined by 
	\[
		\Gamma(S,X)(T)=\Hom[_S]0{S\times T,X}
	\]
	is representable if $S$ is locally in $\cat{fd}$, and $X$ and $F$ are locally in $\cat{met}$.
\end{Prop}

\begin{proof}
	Let $U_i\subset S$ be open subspaces covering $S$, and $\tau_i\in\tau_X^F(U_i)$. For all $U=U_i$, there $\tau=\tau_i$ is an isomorphism $U\times F\to X_U$ over $U$. 
	
	For any $T$, this gives an isomorphism
	\[
		\GHom0{U,F}(T)=\Hom0{U\times T,F}=\Hom[_U]0{U\times T,U\times F}\to \Hom[_U]0{U\times T,X_U}
	\]
	which is natural in $T$. Thus, we obtain $\Gamma(S,X)$ by gluing the $\Gamma(U,X_U)=\GHom0{U,F}$ along these isomorphisms.  
\end{proof}

\begin{Def}
	If $X/S$ is a bundle and the assumptions of \thmref{Prop}{inner-sections} are satisfied, then we say that $\Gamma(S,X)$ exists, and call this the \Define{superspace of sections of $X$}. In this case, for any open subspace $U\subset S$, $\Gamma(U,X):=\Gamma(U,X_U)$ exists, and $U\mapsto\Gamma(U,X)$ is a sheaf on $S$ with values in $\SMan$. 
\end{Def}

\begin{Par}
	Let $E,F$ be locally convex super-vector spaces over $\knums$. Let $\GPol0{E,F}$ be the subset of all $\GHom0{S(E),F}$ be the subset of all $f$ \scth the maps 
	\[
		(e_1,\dotsc,e_k)\mapsto f(e_1\dotsm e_k):E^k\to F
	\]
	are continuous for any $k\in\nats$. 
	
	To any $f\in\GPol0{E,F}$, we may associate an element $\tilde f\in\GCt[^\infty]0{E,F}$ by 
	\[
		\tilde f(P;x)=\sum_{k=0}^\infty\frac1{k!} f(x^kP)\mathfa P\in S(E)\,,\,x\in E .
	\]
	
	Thus, $\GPol0{E,F}$ is, as a subspace of $\GCt[^\infty]0{E,F}$, a locally convex super-vector space, locally convex if $F$ is (and almost never metrisable). In particular, the set $\GLin0{E,F}$ of all continuous linear maps $E\to F$ can be considered as locally convex super-vector space. Observe that $\GPol0{E,F}_0=\sh P(E,F)$. 

	On the other hand, an even continuous $k$-linear map $f:E_1\times\dotsm\times E_k\to F$ may be considered as an element of $\sh P_k(\prod_jE_j,F)=\GPol[_k]0{\prod_jE_j,F}_0$ by setting
	\[
		f(P_1\otimes\dotsm\otimes P_k)=0 \mathfa P_j\in S^{m_j}(E_j)
	\]
	where at least one $m_j\neq1$, and 
	\[
		f(e_1\otimes\dotsm\otimes e_k)=f(e_1,\dotsc,e_k)\mathfa e_j\in E_j\ .
	\]
	
	The corresponding element of $\Ct[_E^\infty]0{E_0,F}$, where $E=\prod_{j=1}^kE_j$, is then given for $P_j\in S^{m_j}(E_j)$ by 
	\[
		\tilde f(P_1\otimes\dotsm\otimes P_k;e_1,\dotsc,e_k)=
		\begin{cases}
			0 & \min_{j=1}^km_j>1\ ,\\
			f\Parens1{e_1^{1-m_j}P_1,\dotsc,e_k^{1-m_j}P_k} & \text{otherwise.}
		\end{cases}
	\]
	
	In particular, the vector space operations 
	\[
		\cdot:\knums\times E\to E\nd +:E\times E\to E
	\]
	are elements of $\sh P_2(\knums\times E,E)$ resp.~$\sh P_2(E\times E,E)$, and may thus be considered as morphisms of linear supermanifolds $\knums\times E\to E$ and $E\times E\to E$, respectively. This gives meaning to the following proposition. 
 \end{Par}

\begin{Prop}[linear-points]
	Let $E,F$ be metrisable locally convex super-vector spaces and $S\in\Ob\SMan$, locally in $\cat{met}$. 
	
	Fix $k\in\nats$ and set $\GPol[_k]0{E,F}=\GPol0{E,F}\cap\GHom0{S^k(E),F}$. There is a natural bijection between $\GPol[_k]0{E,F}(S)$ and the set of $f\in\Hom0{S\times E^k,F}$ \scth 
	\begin{gather}
		f(s,\dotsc,\lambda\cdot x_j+y_j,\dotsc)=\lambda\cdot f(s,\dotsc,x_j,\dotsc)+f(s,\dotsc,y_j,\dotsc)\in_T F\label{eq:multilin}\\
		f(s,\dotsc,x_i,\dotsc,x_j,\dotsc)=f(s,\dotsc,x_j,\dotsc,x_i,\dotsc)\in_TF\label{eq:symm}
	\end{gather}
	\fa $s\in_TS$, $\lambda\in_T\knums$, $x_1,\dotsc,x_k,y_j\in_T E$, $T\in\Ob\SMan$, and $i,j=1,\dotsc,k$.
\end{Prop}

\begin{proof}
	The map $\phi:\GPol[_k]0{E,F}(S)\to\Hom0{S\times E^k,F}$ is defined for $S=G_W$ by 
	\[
		\phi(f)(P\otimes Q;x,y_1,\dotsc,y_k)=\sum_{n_1=\dotsm=n_k=0}^\infty\frac1{n_1!\dotsm n_k!}[f(P;x)](y_1^{n_1}\dotsm y_k^{n_k}\mu^{k-1}(Q))\ .
	\]
	The map is natural in $S$, as can be checked using Equation \eqref{eq:super-faadibruno}. Thus, the claim only has to be proved in the case of a superdomain $S=G_W$. 
	
	To see that the elements in the range of $\phi$ satisfy the Equations \eqref{eq:multilin} and \eqref{eq:symm}, it suffices to take $S=\GPol[_k]0{E,F}$ and $f=\id$. Then 
	\[
		\phi(f)(P\otimes Q;x,y_1,\dotsc,y_k)=\sum_{n_1+\dotsm+n_k=k-\ell}\frac1{n_1!\dotsm n_k!}T(y_1^{n_1}\dotsm y_k^{n_k}e_1\dotsm e_k)
	\]
	where $(P,x)=(1,T)$ or $(P,x)=(T,x)$, and $Q=e_1\otimes\dotsm\otimes1\otimes\dotsm e_k$ where the number of tensor factor equal to $1$ is $\ell\sle k$, and the other factors are $e_j\in E$. The right hand side, considered as a function of $Q$, $y_j$, is just $T\in\GHom0{S^k(E),F}$, so the equations follow. 
	
	Conversely, assume that $S=G_W$ and we are given a morphism $f:S\times E^k\to F$ satisfying Equations \eqref{eq:multilin} and \eqref{eq:symm}. Define $h\in\Ct[^\infty_G]0{W,\GPol[_k]0{E,F}}$ by 
	\[
		[h(P;x)](e_1\dotsm e_k)=f(P\otimes e_1\otimes\dotsm\otimes e_k;x,0,\dotsc,0)
	\]
	It is clear that $h$ is well defined. 
	
	To see that $f=\phi(h)$, consider for $\lambda\in\knums$ the morphism $d_\lambda:E\to E$ given by the linear map $e\mapsto\lambda\cdot e$. Then $d_\lambda\in\Ct[_E^\infty]0{E_0,F}$ is given by 
	\[
		d_\lambda(e_1\dotsm e_n;x)=\begin{cases}0&n>1\ ,\\ \lambda x^{n-1}e_n & n\sle 1\ ,\end{cases}
	\]
	where $e_0=1$. Thus, by applying the definitions, one obtains \fa $n\sge0$, $m\sge1$, 
	\[
		d_{\lambda(n)}(e_1\dotsm e_m;x)=\delta_{n+1,m}\lambda^me_1\dotsm e_m
	\]
	
	We have $\lambda_1\dotsm\lambda_k\cdot f=f\circ(\id\times d_{\lambda_1}\times\dotsm\times d_{\lambda_k})$ by the assumption on $f$, and this gives the equation
	\begin{multline*}
		f(P\otimes Q_1\otimes\dotsm\otimes Q_k;x,y_1,\dotsc,y_k)\\=\lambda_1^{m_1-1}\dotsm\lambda_k^{m_k-1}f(P\otimes Q_1\otimes\dotsm\otimes Q_k;x,\lambda_1y_1,\dotsc,\lambda_ky_k)
	\end{multline*}
	\fa $\lambda_j\neq0$, $P\in S(G)$, $Q_j\in S^{m_j}(E)$, $x\in W$ and $y_j\in E_0$. Hence
	\[
		f(P\otimes Q_1\otimes\dotsm\otimes Q_k;x,y_1,\dotsc,y_k)=0\mathtxt{whenever}\min\nolimits_{j=1}^km_j>1\ .
	\]
	For $m_1=\dotsm=m_k=1$, we get 
	\[
		f(P\otimes Q_1\otimes\dotsm\otimes Q_k;x,y_1,\dotsc,y_k)=f(P\otimes Q_1\otimes\dotsm\otimes Q_k;x,0,\dotsc,0)\ .
	\]
	
	Finally, using the above equation for $m_j=1$ (some $j$), we obtain, using the abbreviation $Q=Q_1\otimes\dotsm\otimes Q_{j-1}\otimes 1\otimes Q_j\otimes\dotsm Q_k$ (where $m_i\sle1$ for $i\neq j$), 
	\begin{align*}
		f(P\otimes Q;x,y_1,\dotsc,y_k)&=\frac d{d\lambda}\Big|_{\lambda=0}f(P\otimes Q;x,y_1,\dotsc,(1+\lambda)y_j,\dotsc,y_k)\\
		&=f(P\otimes Q(1^{j-1}\otimes y_j\otimes 1^{k-j});x,y_1,\dotsc,y_k)\ .
	\end{align*}
	Hence, for $m_i\sle 1$ (all $i$) and $Q_i=1$ at some positions $i=i_1<\dotsm<i_\ell$, 
	\begin{align*}
		f(P\otimes Q_1\dotsm 1\dotsm Q_k;x,y_1,\dotsc,y_k)&=f(P\otimes Q_1\dotsm y_{i_1}\dotsm Q_k;x,y_1,\dotsc,y_k)\\
		&=\dotsm\\
		&=f(P\otimes Q_1\dotsm y_{i_1}\dotsm y_{i_\ell}\dotsm Q_k;x,y_1,\dotsc,y_k)\\
		&=f(P\otimes Q_1\dotsm y_{i_1}\dotsm y_{i_\ell}\dotsm Q_k;x,0,\dotsc,0)
	\end{align*}
	This readily entails that $\phi(h)=f$, and hence, the claim. 
\end{proof}

\begin{Def}[def-vectbun]
	Let $E,F$ be locally convex super-vector spaces and $S\in\Ob\SMan$. A morphism $f\in\Hom[_S]0{S\times E,S\times F}=\Hom0{S\times E,F}$ over $S$ is called \Define{linear (over $S$)} if Equation \eqref{eq:multilin} is fulfilled for $k=1$. 
	More generally, if Equations \eqref{eq:multilin} and \eqref{eq:symm} are fulfilled for $f$, where $k$ is now arbitrary, we say that $f$ is a \Define{homogeneous polynomial morphism of degree $k$}.

	If $X/S$ is a bundle whose fibres are linear supermanifolds, and such that there exists an open cover $(U_i)$ of $S$ and trivialisations $\tau_i$ over $U_i$ \scth $\tau_{ji}=\tau_j^{-1}\circ\tau_i$ is linear over $U_{ij}$, for any $i,j$, then $X/S$ is called a \Define{vector bundle}. A morphism $f:X/S\to Y/S$ of such vector bundles is a \Define{morphism of vector bundles} if after trivialisation, it becomes linear over $S$. 
	
	More generally, if $X/S$ and $Y/T$ are vector bundles and $f:X\to Y$ is a morphism over $\vphi:S\to T$, then $f$ is called a \Define{morphism of vector bundles} if it induces a morphism of vector bundles $X\to S\times_TY$ over $S$. 
\end{Def}

\begin{Prop}
	Let $X/S$ be a vector bundle where $S$ is locally in $\cat{fd}$, and $X$ and the fibre $F$ are locally in $\cat{met}$. Then $\Gamma(-,X)$ is a sheaf of locally convex super-vector spaces on $S$ which are metrisable if $S_0$ is. 
\end{Prop}

\begin{proof}
	On a trivialising open set, this is obvious. In general, we have a Hausdorff quotient of a product of locally convex spaces, so the statement is immediate. In case $S_0$ is metrisable, we have a quotient of a metrisable Hausdorff locally convex space, which is itself second countable and therefore also metrisable. 
\end{proof}

We take note of the following simple fact.

\begin{Lem}
	Let $V/S$ be a vector bundle and $s\in S_0$. Then the fibre $V_s$ is a locally convex super-vector space. 
\end{Lem}

\subsection{Weil bundles} 

As is known \cite{weil,kolar-michor-slovak,koszul-near,fioresi-local}, the formalism of Weil bundles forms a convient setting for the uniform discussion of (iterated) tangent and jet bundles. In the super context, the (higher) odd tangent bundles \cite{gorms-worms} also naturally fall into this framework. As we will show, the Weil bundles can be introduced easily by considering the spectra of Weil superalgebras. These can naturally be realised within the category $\SSp$, and the Weil bundles may then be defined as inner homs. Although this approach may appear counterintuitive at first sight, it is in fact very natural, as will be  illustrated below by our discussion of loop supergroups. 

\begin{Def}
	A \Define{Weil superalgebra} is a finite-dimensional, super-commutative, and local superalgebra over $\reals$. Since $A_1$ consists of non-units, $A_1\subset\ger m$ where $\ger m$ is the maximal ideal, so $R/\ger m=\reals$ and $A=\reals\oplus\ger m$. Since $A$ is Noetherian, there is, by Krull's Intersection Theorem, some $k\in\nats$ \scth $\ger m^{k+1}=0$, \ie $\ger m$ is nilpotent. The minimal such $k$ is called the \Define{height} of $A$. Similarly, $p|q=\dim\ger m/\ger m^2$ is called the \Define{width} of $A$. Any Weil superalgebra of height $k$ and width $p|q$ is isomorphic to a quotient of $\reals[x_1,\dotsc,x_p|\xi_1,\dotsc,\xi_q]/(x_1^{k+1},\dotsc,x_p^{k+1})$ where $\reals[x_1,\dotsc,x_p|\xi_1,\dotsc,\xi_q]=\reals[x_1,\dotsc,x_p]\otimes\bigwedge(\xi_1,\dotsc,\xi_q)$ denotes the superpolynomial algebra in $p|q$ homogeneous generators.
\end{Def}

\begin{Par}
	Let $A$ be a Weil superalgebra with its unique Hausdorff locally convex topology. We define $\Spec A\in\SSp$ as follows. The underlying topological space $(\Spec A)_0$ is the point $*$. For $E_U\in\SPol$, let $\sh O_{\Spec A}(*,E_U)=U\times(\ger m\otimes E)_0$, with the inclusion $U\subset U\times(\ger m\otimes E)_0$ given by $u\mapsto(u,0)$. 
	
	Next, let $f:E_U\to F_V$ is a morphism in $\SPol$. To define $\sh O_{\Spec A}(*,f)$, observe that $*=0$ is open in $0$, and that $U\times(\ger m\otimes E)_0$ is open in $(A\otimes E)_0$. (Here, note that any element on $A\otimes E$ can be expressed uniquely as $q\otimes x+y$ where $x\in E$ and $y\in\ger m\otimes E$; so $U\times (\ger m\otimes E)_0$ identifies with the open subset of $(A\otimes E)_0$ of all $t=1\otimes x+y\in(A\otimes E)_0$ \scth $x\in U$.) Moreover,
	\[
		\Ct[^\infty_{0,A\otimes E}]0{*,U\times(\ger m\otimes E)_0}=U\times(\ger m\otimes E)_0\ .
	\]
	
	The map $\id_A\otimes f:S(A\otimes E)\to S(F)$ is in $\sh P(A\otimes E,A\otimes F)$, and hence in $\Ct[^\infty_{A\otimes E}]0{(A\otimes E)_0,A\otimes F}$. Since the underlying map sends $a\otimes e$ to $a\otimes f(e)$, it is clear that $(\id_A\otimes f)\circ(-)$ makes sense as a map
	\[
		U\times(\ger m\otimes E)_0\to V\times(\ger m\otimes F)_0\ .
	\]
	Hence, we may define $\sh O_A(*,f)=(\id_A\otimes f)\circ(-)$. 
\end{Par}

\begin{Rem}
	An example of a Weil superalgebra is $A=\bigwedge(\reals^q)^*$. In this case, $\Spec A=\reals^{0|q}$, and the above construction coincides with the embedding of this finite-dimensional supermanifold into the category $\SMan$. 
\end{Rem}

\begin{Prop}[weilprod-exists]
	Let $X\in\Ob\SSp$ and $A$ be a Weil superalgebra. Then the product $X\times\Spec A$ exists in $\SSp$. 
\end{Prop}

\begin{proof}
	Let $(X\times\Spec A)_0=X$. For any morphism $f:E_U\to F_V$ in $\SPol$, define 
	\[
		\sh O_{X\times\Spec A}(-,E_U)=\sh O_X(-,(U\times(\ger m\otimes E)_0,A\otimes E))
	\]
	and $\sh O_{X\times\Spec A}(-,f)=\,\sh O_X(-,\id_A\otimes f)$.
	
	The second projection $p_2:X\times \Spec A\to\Spec A$ is defined on the level of the structure functor by the inclusions 
	\[
		U\times(\ger m\otimes E)_0\subset\sh O_X(-,(U\times(\ger m\otimes E)_0,A\otimes E))\ . 
	\]
	
	The first projection $p_1:X\times\Spec A\to X$ is defined by applying the functor $\sh O_X$ to the obvious superpolynomial map $E_U\to(U\times(\ger m\otimes E)_0,A\otimes E)$. 
	
	The non-trivial point we need to prove in order to establish that the superspace $X\times\Spec A$, as it is defined, really is a product in the category $\SSp$, is that given morphisms $f_1:Y\to X$ and $f_2:Y\to\Spec A$, there exists a unique morphism $f:Y\to X\times\Spec A$ \scth $p_j\circ f=f_j$, $j=1,2$. 
	
	To see the existence, let $\phi=(f_2)^*_{\reals,\reals^{1|1}}:A\to\sh O_Y(Y_0,(\reals,\reals^{1|1}))$. By naturality of $f_2^*$, $\phi$ is an even unital algebra morphism, where the superalgebra structure on $\sh O_Y(Y_0,(\reals,\reals^{1|1}))$ is induced by that on $\reals^{1|1}$. For any $e\in E$, consider the even linear map $p_e:\reals^{1|1}\to E$ given by $p_e(x,y)=xe_0+ye_1$. Applying naturality to this, we see 
	\[
		(f_2)^*_{E_U}\Parens1{\textstyle\sum_ja_j\otimes e_j}=\textstyle\sum_jp_{e_j}\circ \phi(a_j)
	\]
	
	For $E_U\in\Ob\SPol$, $\sh O_X(-,(E_0,E))$ is a sheaf of super-vector spaces, with structure obtained by the functoriality of $\sh O_X$ from that of $E$. Moreover, the product sheaf $\sh O_X(-,E_U)\times(\ger m\otimes\sh O_X(-,(E,E^{1|1})))_0$ (where $E^{1|1}=E\oplus\Pi E$) may naturally be considered as a subsheaf of $(A\otimes\sh O_X(-,(E,E^{1|1})))_0$. We will now define sheaf morphisms 
	\[
		\xymatrix{%
			\sh O_X(-,(U\times(\ger m\otimes E)_0,A\otimes E))\ar@<+0.5ex>[r]^-{\psi}&\sh O_X(-,E_U)\times(\ger m\otimes \sh O_X(-,(E,E^{1|1})))_0\ar@<+0.5ex>[l]^-{\vphi}}
	\]
	which are mutually inverse to each other. 
	
	To that end, let $a_0=1,a_1,\dotsc,a_n$ be a homogeneous basis of $A$ where $a_j\in\ger m$ for $j>0$, and $\alpha_0,\alpha_1,\dotsc,\alpha_n$ the dual basis of $A^*$. We obtain even linear maps 
	\[
		\alpha_j\otimes\id:A\otimes E\to E^{1|1}\nd a_j\otimes\id:\Pi^jE\to A\otimes E
	\]
	which may be considered as a morphisms in $\SPol$. Define  
	\[
		\psi=\textstyle\sum_ja_j\otimes\sh O_X(-,\alpha_j\otimes\id)\nd\vphi=\textstyle\sum_j\alpha_j\otimes\sh O_X(-,a_j\otimes\id)\ .
	\]
	
	 By naturality, these sheaf morphisms are well defined, and it is routine to check that they are mutually inverse. One may also see that they are independent of the choice of basis, so that there is a natural isomorphism 
	 \[
	 	\sh O_X(-,(U\times(\ger m\otimes E)_0,A\otimes E))\cong\sh O_X(-,E_U)\times(\ger m\otimes\sh O_X(-,(E,E^{1|1})))_0\ .
	\]
	
	It is now clear how to define $f:Y\to X\times\Spec A$ \scth $p_j\circ f=f_j$: Let $f_0=f_{1,0}$ and define 
	\[
		f^*_{E_U}:\sh O_X(-,(U\times(\ger m\otimes E)_0,A\otimes E))\to\sh O_Y(-,E_U)
	\]
	by $f^*_{E_U}\Parens1{\textstyle\sum_j a_j\otimes h_j}=\textstyle\sum_j\phi(a_j)\cdot(f_1)^*_{E_U}(h_j)$. Here, the action of $\sh O_Y(-,\reals^{1|1})$ on $\sh O_Y(-,(E_0,E))$ is induced \via the functoriality of $\sh O_Y$ from that of $\reals^{1|1}$ on $E$. Somewhat tedious verifications show that this defines a morphism $f$ with the required properties. 
	
	The uniqueness of $f$ can be seen by applying similar ideas to the even linear isomorphism $A\cong(\reals^{1|1}\otimes A)_0$. Hence follows the assertion.
\end{proof}

\begin{Lem}[spec-hom]
	For any $X\in\Ob\SMan$ and any superdomain $F_V$, we have a natural isomorphism
	\[
		\Hom0{X\times\Spec A,F_V}=\sh O_{X\times\Spec A}(X_0,F_V)
	\]
\end{Lem}

\begin{proof}
	There is a natural map from the left to the right hand side, given by 
	\[
		f\mapsto f^*(\id_{F_V})
	\]
	To show that it is an isomorphism, we may assume that $X=E_U$, the general case following by naturality. 
	
	Let $f\in\sh O_{X\otimes\Spec A}(X_0,F_V)=\Ct[^\infty]0{E_U,(A\otimes F)_{U\times(\ger m\otimes E)_0}}$. The point is to define $f^*(h)\in\Ct[^\infty]0{E_U,(A\otimes G)_{W\times(\ger m\otimes G)_0}}$ for any $h\in\Ct[^\infty]0{F_V,G_W}$. We will set 
	\[
		f^*(h)=h_A\circ f\ ,
	\]
	once having defined $h_A\in\Ct[^\infty]0{(A\times F)_{V\times(\ger m\otimes F)_0},(A\otimes G)_{W\times(\ger m\otimes G)_0}}$. This is a version of the usual `Taylor expansion in nilpotent directions'. 
	
	First, we define $h(a\otimes P;x)=a\otimes h(P;x)$ for $a\in A$, $P\in S(F)$, $x\in V$, and extend this $A$-linearly. Then we set 
	\begin{equation}\label{eq:weilfun-local}
		h_A(P;x+y)=\sum_{j=0}^\infty\frac1{j!}h_A(y^jP;x)
	\end{equation}
	\fa $P\in S(A\otimes F)$, $x\in V$, $y\in(\ger m\otimes F)_0$. A number of routine checks shows that this proves the assertion. 
\end{proof}

\begin{Cor}[weil-innerhom-domain]
	Let $E_U$ be a superdomain and $A$ a Weil superalgebra. Then $\GHom0{\Spec A,E_U}$ exists in $\SMan$, and is given by $(A\otimes E)_{U\times (\ger m\otimes E)_0}=E_U\times \ger m\otimes E$. 
\end{Cor}

\begin{proof}
	For any supermanifold $X$, we have natural isomorphisms
	\begin{align*}
		\Hom0{X\times\Spec A,E_U}&=\sh O_{X\times\Spec A}(X_0,E_U)\\
		&=\sh O_X(X_0,(U\times(\ger m\otimes E)_0,A\otimes E))\\
		&=\Hom0{X,(A\otimes E)_{U\times(\ger m\otimes E)_0}}\ ,
	\end{align*}
	in view of \thmref{Lem}{spec-hom}, \thmref{Cor}{super-mor}, \thmref{Prop}{weilprod-exists}, and its proof. 
\end{proof}

In order to extend this local existence result, we introduce some natural constructions with $\GHom0{\Spec A,-}$.

\begin{Par}
	Let $A$ be a Weil superalgebra and $X,Y,Z$ be supermanifolds. We write $T^AX=\GHom0{\Spec A,X}$ whenever it exists (as a supermanifold). 
	
	Assume that $T^AX$ and $T^AY$ exist, and that $\vphi:X\to Y$ is a morphism. We define $T^A\vphi:T^AX\to T^AY$ to be the morphism given on $S$-points by:
	\[
		T^A\vphi(s)=\vphi(s)\in_{S\times\Spec A}Y\ .
	\]
	Here, observe that $\GHom0{\Spec A, X}(S)=X(S\times\Spec A)$. This construction is obviously functorial, \emph{i.e.}~$T^A\id_X=\id_{T^AX}$, and if $T^AZ$ exists and $\psi:Y\to Z$ is a morphism, then $T^A(\psi\circ\vphi)=T^A\psi\circ T^A\vphi$. 
	
	If $B$ is another Weil superalgebra, and $f:A\to B$ is an even unital algebra morphism, then we may define a natural transformation $T^f:T^A\to T^B$ (defined where both functors are defined), as follows. To $f$ there is associated a morphism $f^*:\Spec B\to\Spec A$: Its underlying map is the identity, and on the level of structure functors, 
	\[
		\sh O_{\Spec A}(*,E_U)=U\times(\ger m_A\otimes E)_0\to\sh O_{\Spec B}(*,E_U)= U\times(\ger m_B\otimes E)_0
	\] 
	is given by $f\otimes\id_E$ (observe that any morphism of Weil superalgebras is automatically local). Now $T^f_X:T^AX\to T^BX$ is defined on $S$-points by 
	\[
		T^f_X(s)=s\circ(\id_S\times f^*)\in X(S\times\Spec B)=(T^BX)(S)
	\]
	\fa $s\in(T^AX)(S)=X(S\times\Spec A)$.
\end{Par}

\begin{Lem}
	For any even unital algebra morphism $f:A\to B$ of Weil superalgebras, $T^f:T^A\to T^B$ is, on its domain of definition, a natural transformation.
\end{Lem}

\begin{proof}
	Let $\vphi:X\to Y$ be a morphism of supermanifolds \scth $T^CZ$ exists for $C=A,B$ and $Z=X,Y$. Then 
	\begin{align*}
		(T_Y^f\circ T^A\vphi)(s)&=T^f_Y(\vphi(s))=\vphi(s)\circ(\id_S\times f^*)\\
		&=\vphi(s\circ(\id_S\times f^*))=T^B\vphi(s\circ(\id_S\times f^*))=(T^B\vphi\circ T^f_X)(s)
	\end{align*}
	\fa $s\in T^AX(S)=X(S\times\Spec A)$. Hence, $T^f$ is a natural transformation.
\end{proof}

\begin{Par}
	Let $X,Y$ be supermanifolds and $A$ a Weil superalgebra \scth $T^AX$, $T^AY$ exist. Clearly, $\reals$ is a Weil superalgebra, $\Spec\reals=*$, and $\GHom0{\Spec\reals,Z}=Z$ for any superspace $Z$, and $T^\reals\vphi=\vphi$ for any morphism $\vphi$. 
	
	There are two natural morphisms $\eta:\reals\to A$ and $\eps:A\to\reals$, determined uniquely by the requirement of unitality. We let $0_{T^AX}=T^\eta_X$ and $\pi_{T^AX}=T^\eps X$. 

	By the lemma, there is for any morphism $\vphi:X\to Y$ a commutative diagram 
	\[
		\xymatrix@C+2ex{%
			X\ar[r]^-{0_{T^AX}}\ar[d]_-{\vphi}&{\GHom0{\Spec A,X}}\ar[r]^-{\pi_{T^AX}}\ar[d]^-{T^A\vphi}&X\ar[d]^-{\vphi}\\
			Y\ar[r]^-{0_{T^AY}}&{\GHom0{\Spec A,Y}}\ar[r]^-{\pi_{T^AY}}&Y
		}
	\]
	
	In particular, $T^A\vphi$ is a relative morphism over $\vphi$, where we consider $\pi_{T^AX}$ as the relative superspace $T^AX/X$.
\end{Par}

\begin{Prop}[weil-innerhom]
	For any Weil superalgebra $A$ and any supermanifold $X$, the inner hom $T^AX=\GHom0{\Spec A,X}$ exists as a supermanifold. With the morphism $\pi_{T^AX}$, It has the structure of a fibre bundle over $X$. If $E_U$ is a superdomain which is an open subspace of $X$, then $T^AX$, restricted to $E_U$, has fibre $\ger m\otimes E$. 
\end{Prop}

\begin{proof}
	Let $(U_i)$ be an open cover of $X$ by open subspaces isomorphic to superdomains, where $U_i$ is isomorphic to an open subspace of the locally convex super-vector space $E_i$ (say). Then $T^AU_i$ and equals $U_i\times E_i\otimes \ger m$, by \thmref{Cor}{weil-innerhom-domain}. Applying the above construction, we obtain isomorphisms $(T^AU_j)_{U_{ij}}\to (T^AU_i)_{U_{ij}}$ over $U_{ij}$. Then the pasting lemma for bundles implies our claim. 
\end{proof}

\begin{Def}
	Let $A$ be a Weil superalgebra. The functor $T^A$ from supermanifolds to bundles is called the \Define{$A$-Weil functor}. In particular, for $A=\reals[\eps]/(\eps^2)$, it is called the \Define{tangent functor} and denoted by $T$. 
\end{Def}

\begin{Prop}[weilfunc-basic]
	Let $A,B$ be Weil superalgebras. The Weil functors $T^A,T^B$ enjoy the following properties:
	\begin{enumerate}
		\item $T^A$ commutes with finite products.
		\item If $U$ is an open subspace of the supermanifold $X$, then $T^AU$ is an open subspace of $T^AX$.
		\item There is a natural isomorphism $T^{A\otimes B}\cong T^AT^B$; in particular, there is a natural isomorphism $T^A\circ T^B\cong T^B\circ T^A$. 
	\end{enumerate}
\end{Prop}

\begin{proof}
	(1). For the empty product, we have 
	\[
		\Hom0{S,T_A(*)}=\Hom0{S\times\Spec A,*}=*\ ,
	\]
	so $T_A(*)=*$. In the case of binary products, this may by definition be checked locally. But then it is obvious by \thmref{Cor}{weil-innerhom-domain}. 
	
	(2). This is clear by construction.
	
	(3). We have natural isomorphisms 
	\begin{align*}
		\Hom0{S,T^{A\otimes B}X}&=\Hom0{S\times\Spec A\otimes B,X}\\
		&=\Hom0{S\times\Spec A,T^BX}=\Hom0{S,T^A(T^BX)}
	\end{align*}
	for any supermanifold $S$. The second equality holds, since 
	\begin{align*}
		\Hom0{X\times\Spec A\otimes B,E_U}&=\sh O_{X\times\Spec A\otimes B}(-,E_U)\\
		&=\sh O_X(-,(U\times (\ger m_{A\otimes B}\otimes E),A\otimes B\otimes E))\\
		&=\sh O_{X\times\Spec A}(-,(U\times(\ger m_B\otimes E)_0,B\otimes E))\\
		&=\Hom0{X\times\Spec A,T^BE_U}
	\end{align*}
	for any superdomain $E_U$. This proves the assertion.
\end{proof}

\begin{Prop}
	For any $X\in\Ob\SMan$, $TX$ is a vector bundle. More generally, if $A$ is a Weil superalgebra \scth $\ger m^2=0$, then $T^AX$ is a vector bundle.
\end{Prop}

\begin{proof}
	It is sufficient to prove that $T^Af$ is linear for any morphism $f:E_U\to F_V$ of superdomains. But by construction, $T^Af=f_A$, the latter being defined in Equation \eqref{eq:weilfun-local}. We find
	\[
		T^Af(P+Q;x+y)=f(P;x)+f(Q;x)+f(yP;x)
	\]
	\fa $P\in S(E)$, $Q\in S(\ger m\otimes E)$, $x\in U$, $y\in(\ger m\otimes E)_0$. According to the proof of \thmref{Prop}{linear-points}, this shows that $T^Af$ is linear.
\end{proof}

\begin{Rem}
	The converse statement is also true: If $A$ is a Weil superalgebra \scth $T^AX$ is a vector bundle \fa $X\in\Ob\SMan$, then $\ger m^2=0$. 
\end{Rem}

We end this section by a brief discussion of vector fields. 

\begin{Def}
	Let $X\in\Ob\SMan$ be locally in $\cat{fd}$. The super-vector space $\Gamma(X,TX)$ is denoted $\ger X(X)$. Its elements are called \Define{vector fields}. 
\end{Def}

\begin{Par}
	Let $X\in\Ob\SMan$ be locally in $\cat{fd}$, and $E$ locally convex super-vector space. We will now define a morphism $\ger X(X)\times\GHom0{X,E}\to\GHom0{X,E}$ which is bilinear and natural in $X$ and $E$. Recall that $\GHom0{X,E}=\sh O_X(X_0,E^{1|1})$, so this morphism we be a generalisation of the usual action of vector fields on functions. 
	
	Let $S\in\Ob\SMan$ and $x\in_S\ger X(X)$, $f\in_S\GHom0{X,E}$, so that we may consider $x\in\Hom[_S]0{S\times X,S\times TX}$ and $f\in\Hom[_S]0{S\times X,S\times E}$. 
	
	Define $T_Sf\in\Hom[_S]0{S\times TX,S\times TE}$ by 
	\[
		T_Sf=(\pi_{TS}\times\id_{TE})\circ Tf\circ(0_{TS}\times\id_{TX})\ ,
	\]
	\ie $T_Sf$ is the partial tangent in the $X$ variable. Then $T_Sf$ is a relative morphism over $\id_S\times f$, and $TE=E\times E$, so that we may consider
	\[
		T_Sf\in\Hom[_{S\times X}]0{S\times TX,X\times_ETE}=\Hom[_S]0{S\times TX,E}\ .
	\]
		
 	Now, we may define $x(f)\in_S\GHom0{X,E}$ by $x(f)=T_Sf\circ x$. Since $T_Sf$ is a vector bundle morphism, and the vector bundle structure of $TE$ comes from the linear structure on $E$, it is easy to check that this indeed defines a bilinear morphism. The naturality is obvious by construction. 
\end{Par}

\begin{Def}
	If $\ger g$ is a Lie superalgebra in $\SMan$, a \Define{smooth representation of $\ger g$ on $V$} is an even continuous bilinear map $\alpha:\ger g\times V\to V$ \scth 
	\[
		\alpha(x,\alpha(y,v))-(-1)^{\Abs0x\Abs0y}\alpha(y,\alpha(x,v))=\alpha([x,y],v)
	\]
	for all homogeneous $x,y\in\ger g$, and all $v\in V$. 
\end{Def}

\begin{Prop}
	Let $X\in\Ob\SMan$ be locally in $\cat{fd}$, and $E$ locally convex super-vector space. The canonical morphism 
	\[
		\ger X(X)\times\GHom0{X,E}\to\GHom0{X,E}
	\]
	defines a representation of the Lie superalgebra $\ger X(X)$ on the locally convex super-vector space $\GHom0{X,E}$.
\end{Prop}

We omit the proof, since we will not use the result below. 

\subsection{Lie supergroups and their representations}

\begin{Def}
	A \Define{Lie supergroup} is a group object in $\SMan$.
\end{Def}

\begin{Par}
	If $f\in\Ct[^\infty_{E,F}]0{U,V}$ where $E_U,F_V\in\SDom_{\sh C^\infty}$, then we may consider $f|_{S(E_0)}=f_0\in\Ct[^\infty_{E_0,F}]0{U,V}=\Ct[^\infty]0{U,V}$. The map $f\mapsto f_0$ is compatible with composition, \ie$(g\circ f)_0=g_0\circ f_0$ for $g\in\Ct[_{F,G}^\infty]0{V,W}$. 
	
	Hence, if $X\in\SSp^\reals$ has an open cover by superdomains, there is canonical superspace $X_0\in\SSp^\reals$ on the topological space $X_0$, and a canonical morphism $j_0:X_0\to X$ \scth $j_0^*(f)=f_0$ whenever $U\subset X_0$ is open \scth $X_U$ is a superdomain, and $f\in\sh O_X(U,F_V)$. 
\end{Par}

\begin{Lem}
	The correspondence $X\mapsto X_0$ defines a functor $[-]_0$, and the morphisms  $j_0:X_0\to X$ define a natural transformation $[-]_0\to\id$. On supermanifolds, $[-]_0$ commutes with finite products. In particular, if $G$ is a supergroup, then $G_0$ is a Lie group in $\SMan$. 
\end{Lem}

\begin{Def}
	Let $X\in\Ob\SMan$ and $G$ a Lie supergroup. A \Define{left action of $G$ on $V$} is a morphism $\alpha:G\times V\to V$ which satisfies the equations 
	\[
		\alpha(g_1,\alpha(g_2,x))=\alpha(g_1g_2,x)\nd\alpha(1,x)=x
	\]
	\fa $g_1,g_2\in_TG$, $x\in_TX$, and $T\in\Ob\SMan$. One also says that $X$ is a (left) \Define{$G$-space}, and writes $g.x$ for $\alpha(g,x)$. Similarly, one defines right actions (or, equivalently, right $G$-spaces).
	
	A morphism $f:X\to Y$ of $G$-spaces is \Define{$G$-equivariant} if 
	\[
		f(\alpha_X(g,v))=\alpha_Y(g,f(x))\mathfa g\in_SG\,,\,x\in_SX\ ,
	\]
	and all $S\in\Ob\SMan$. 
	
	Let $V$ be a vector bundle over a $G$-space $X$. Then $V$ is a \Define{$G$-equivariant vector bundle} if there is given an action 
	\[
		\alpha_V:G\times V\to V
	\]
	which is linear over $G\times X$, \scth the bundle projection $\pi:V\to X$ is equivariant. 
	
	In the case that $V$ is a trivial vector bundle over $*$, we say that $\alpha_V$ is a \Define{representation of $G$}. To emphasise that morphisms are smooth, we will sometimes call this a \Define{smooth representation}.
\end{Def}

\begin{Def}
	Let $G$ be a Lie supergroup. The locally convex super-vector space $T_1G=(TG)_1$ is denoted by $\ger g$ and called the \Define{Lie superalgebra} of $G$. 
\end{Def}

\begin{Prop}
	Let $G$ be a Lie supergroup. The tangent multiplication on $TG$ induces the structure of a $G$-equivariant vector bundle on $TG$, and as such, $TG$ and $G\times\ger g$ are naturally isomorphic. 
\end{Prop}

\begin{proof}
	The multiplication $m:G\times G\to G$ defines a left action of $G$ on itself. Clearly, $Tm$ defines a left action of $TG$ on itself which is linear over $TG\times G$. On the other hand, $0_{TG}:G\to TG$ is a morphism of Lie supergroups, so that we obtain a left action $\alpha_{TG}:G\times TG\to TG$ with respect to which $TG$ is indeed a $G$-equivariant vector bundle over $G$. 
	
	As $\ger g=T_1G=TG\times_G1$, there is a natural monomorphism $\pi_1:\ger g\to TG$. Consider the composite
	\[
		\xymatrix@C+3ex{%
			\vphi:G\times\ger g\ar[r]^-{\id\times\pi_1}&G\times TG\ar[r]^-{\alpha_{TG}}&TG
		}
	\]
	This is manifestly a morphism of $G$-equivariant vector bundles. 
	
	To define a morphism $f:TG\to\ger g$, consider $f_2:TG\to 1$ and $f_1:TG\to TG$, given by the composite
	\[
		\xymatrix@C+3ex{%
			f_1:TG\ar[r]^-{(\pi_{TG},\id)}&G\times TG\ar[r]^-{i\times\id}&G\times TG\ar[r]^-{\alpha_{TG}}&TG
		}	
	\]
	Then $\pi\circ f_1=1$, so there exists a unique $f:TG\to G$ \scth $\pi_1\circ f=f_1$. 
	
	Now, we define $\psi:TG\to G\times\ger g$ by $\psi=(\pi_{TG},f)$. Then 
	\[
		\vphi\circ\psi=\alpha_{TG}\circ\Parens1{\pi_{TG},\alpha_{TG}\circ(i\times\id)\circ(\pi_{TG},\id)}\ ,
	\]
	so that \fa $x\in_STG$, 
	\[
		\vphi(\psi(x))=\pi_{TG}(x).(\pi_{TG}(x)^{-1}.x)=x\ .
	\]
	
	Conversely, 
	\begin{align*}
		(\id\times\pi_1)&\circ\psi\circ\vphi=(\pi_{TG},f_1)\circ\alpha_{TG}\circ(\id\times\pi_1)\\
		&=\Parens1{m\circ(\id\times\pi_{TG})\circ(\id\times\pi_1),\alpha_{TG}\circ(i\times\id)\times(\pi_{TG},\id)\circ\alpha_{TG}\circ(\id\times\pi_1)}
	\end{align*}
	Because $\pi_{TG}\circ\pi_1=1$, the first argument evaluates to $\id_G$. In the second argument, we compute for $(g,x)\in_S G\times\ger g$, 
	\[
		\pi_{TG}(g.\pi_1(x))^{-1}.(g.\pi_1(x))=g^{-1}.(g.\pi_1(x))=\pi_1(x)\ .
	\]
	Hence, $(\id\times\pi_1)\circ\psi\circ\vphi=\id\times\pi_1$, and since $\id\times\pi_1$ is a monomorphism, we conclude that $\vphi$ is indeed an isomorphism. 
\end{proof}

\begin{Par}
	Let $G$ be a Lie supergroup. \emph{Via} the isomorphism $TG\cong G\times\ger g$, the latter is turned into a Lie supergroup. We define a map $\Ad:G\times\ger g\to\ger g$ by the equation
	\[
		(g,0)(0,x)(g,0)^{-1}=(1,\Ad(g,x))
	\]
	\fa $S$-points $(g,x)\in_SG\times\ger g$. In other words, $\Ad$ is given by the conjugation in $TG\cong G\times\ger g$. We write $\Ad(g)(x)=\Ad(g,x)$ and call $\Ad$ the \Define{adjoint action}.
	
	With this notation, the Lie supergroup structure of $G\times\ger g$ is given on $S$-points $(g,x),(h,y)\in_SG\times\ger g$ by 
	\[
		(g,x)(h,y)=(gh,\Ad(h^{-1})(x)+y)\ ,\ (g,x)^{-1}=(g^{-1},-\Ad(g)(x))\ ,
	\]
	and $1_{G\times\ger g}=(1,0)$. Moreover, $\Ad$ is linear over $G$, so it is a representation on $\ger g$. 
	
	Indeed, $0_{TG}:G\to G\times\ger g$ is a morphism of Lie supergroups, and $Tm$ induces a linear map $\bullet:\ger g\times\ger g\to\ger g$ \scth $0\in\ger g_0$ is neutral. Then 
	\[
		x\bullet y=(x+0)\bullet(0+y)=(x\bullet 0)+(0\bullet y)=x+y\ .
	\]
	In other words, the canonical map $\pi_1:\ger g\to TG$ (which under the isomorphism $TG\cong G\times\ger g$ corresponds to $(1,\id):\ger g\to G\times\ger g$) is a morphism of Lie supergroups, where the supergroup structure on $\ger g$ is induced by the vector space addition. 
\end{Par}

\begin{Def}
	Let $\alpha:G\times V\to V$ be a smooth $G$-representation. We define $d\alpha:\ger g\times V\to V$ as the composite
	\[
		\xymatrix@C+3ex{%
			\ger g\times V\ar[r]^-{\pi_1\times\pi_0}&TG\times TV\ar[r]^-{T\alpha}&TV\ar[r]^-{\pi_{TV}}&V
		}
	\]
	where $\pi_0:V\cong(TV)_0\to TV$ is the canonical map. Under the identification $TV=V\times V$, $\pi_0$ corresponds to $(0,\id):V\to V\times V$. Thus, $d\alpha$ is given on $S$-points $(x,v)\in_S\ger g\times V$ by $(0,d\alpha(x,v))=T\alpha(1,x,0,v)$, so that $d\alpha$ is bilinear. 
	
	In particular, we have a bilinear map $d\Ad:\ger g\times\ger g\to\ger g$, which we denote by $[\cdot,\cdot]$. 
\end{Def}
	
\begin{Par}
	The equation on $S$-points, 
	\[
		\alpha(g,\alpha(h,v))=\alpha(ghg^{-1},\alpha(g,x))
	\]
	gives a corresponding equation on $S\times\Spec A$-points where $A=\reals[\eps]/\eps^2$. 
	
	This implies
	\[
		\alpha(g,d\alpha(x,v))=d\alpha(\Ad(g,x),\alpha(g,v))
	\]
	and
	\[
		d\alpha(x,d\alpha(y,v))=d\alpha([x,y],v)+d\alpha(y,d\alpha(x,v))
	\]
	on $S$-points $g\in_SG$, $x,y\in_S\ger g$, $v\in_SV$. Expanding in Equation \eqref{eq:weilfun-local} with the correct parities gives the equations in the more familiar form
	\[
		\alpha(g,d\alpha(x,v))=d\alpha(\Ad(g,x),\alpha(g,v))
	\]
	and 
	\[
		d\alpha(x,d\alpha(y,v))=d\alpha([x,y],v)+(-1)^{\Abs0x\Abs0y}d\alpha(y,d\alpha(x,v))
	\]
	\fa $g\in G_0$, homogeneous $x,y\in\ger g$, and $v\in V$. In particular, $\ger g$ is a Lie superalgebra. 
\end{Par}

\begin{Def}
	A \Define{supergroup} is given by three data: A Lie superalgebra $\ger g$; a Lie group $G_0$; and a smooth linear action $\Ad:G_0\times\ger g\to\ger g$ by even Lie superalgebra automorphisms, subject to the following conditions: The Lie algebra $\ger g_0$ of $G_0$ is the even part of $\ger g$; the differential $d\Ad$ of $\Ad$ coincides on $\ger g_0\times\ger g$ with the bracket $[\cdot,\cdot]$ of $\ger g$; the adjoint action of $G_0$ coincides on $G_0\times\ger g_0$ with $\Ad$. 
	
	Given a supergroup pair $(\ger g,G_0)$, a \Define{smooth representation} of $(\ger g,G_0)$ is given by a locally convex super-vector space $V$, a smooth representation $\alpha_0:G_0\times V\to V$ of $G_0$, and a smooth representation $d\alpha:\ger g\times V\to V$ of $\ger g$, subject to the following conditions: The differential $d\alpha_0$ of $\alpha_0$ coincides on $\ger g_0\times V$ with $d\alpha$ and 
	\[
		\alpha_0(g,d\alpha(x,v))=d\alpha_0(\Ad(g,x),\alpha_0(g,v))
	\]
	\fa $g\in G_0$, $x\in\ger g$, and $v\in V$. 
\end{Def}
	
We summarise our above considerations as follows.

\begin{Prop}
	Let $G$ be a Lie supergroup. Then $(\ger g,G_0)$, together with the restriction of the adjoint action to $G_0\times\ger g$, is a supergroup pair. Moreover, if $(V,\alpha)$ is a representation of $G$, then $(V,d\alpha,\alpha_0)$ is a smooth representation of $(\ger g,G_0)$, where $\alpha_0=\alpha|_{G_0\times V}$, and $d\alpha$ is the differential of $\alpha$.
\end{Prop}

\begin{Rem}
	For finite-dimensional Lie supergroups, the Lie superalgebra structure on $\ger g$ can be defined along the usual lines, \via left-invariant vector fields. Otherwise, it is not clear how to define the space of vector fields, in the first place. The above derivation circumvents this difficulty. 
\end{Rem}

%
%
%

\begin{Par}
	Let $(\ger g,G_0)$ be a supergroup pair. Consider the canonical super-symmetrisation map $\beta:S(\ger g)\to \Uenv0{\ger g}$, \vq Ref.~\cite{scheunert-lsa}. It is an isomorphism of filtered super-vector spaces, and induces an isomorphism $\Uenv0{\ger g_0}\otimes\bigwedge(\ger g_1)\to\Uenv0{\ger g}$ of filtered $\Uenv0{\ger g_0}$-supermodules. Thus, for any open $V\subset G_0$, and any $E_U\in\Ob\SPol$, we may identify $\Ct[^\infty]0{V\times\ger g_1,E_U}$ with the set $\sh O(V,E_U)$ of all 
	\[
		f=f(u;g)\in\Hom[_{\Uenv0{\ger g_0}}]0{\Uenv0{\ger g},\Ct[^\infty]0{V,E}}
	\]
	\scth $f_0(g)=f(1;g)\in V$ \fa $g\in V$, and \scth \fa $k$, the map 
	\[
		V\times \ger g^k\to E:(g,x_1,\dotsc,x_k)\mapsto f(x_1\dotsm x_k;g)
	\]
	is smooth. 
	
	Following Koszul \cite{koszul}, we introduce on the supermanifold $G_0\times\ger g_1=(G_0,\sh O)$ the structure of a Lie supergroup. Define $1=(1,0):*\to G_0\times\ger g_1$, so that $1^*f=f(1;1)$ for $f\in\sh O(V,E_U)$. Further, set 
	\[
		m^*f(u\otimes v;g,h)=f(\Ad(h^{-1})(u)v;gh)\nd i^*f(u;g)=f(\Ad(g)(S(u));g^{-1})\ .
	\]
	Here, $S:\Uenv0{\ger g}\to\Uenv0{\ger g}$ is defined uniquely by $S(1)=1$, $S(x)=-x$ for $x\in\ger g$, and 
	\[
		S(uv)=(-1)^{\Abs0u\Abs0v}S(vu)\ .
	\]
	
	It is somewhat tedious if straightforward to check the following statement.
\end{Par}

\begin{Prop}
	Endowed with the above structure, $G_0\times\ger g_1$ is a Lie supergroup. 
\end{Prop}

\begin{Par}
	Let $G$ be a finite-dimensional Lie supergroup. To relate it to the supergroup pair $(\ger g,G_0)$, we need to consider the natural action of $\ger g$ on $\sh O_G$. We evade doing this with vector fields, as follows. 
	
	Let $E$ be a metrisable locally convex super-vector space. For $S$-points $g,h\in_SG$, $f\in_S\GHom0{G,E}=\sh O_G(G_0,E^{1|1})$, define
	\[
		(g.f)(h)=f(hg)\in_SE\ .
	\]
	This defines a morphism in 
	\[
		\Hom0{G\times G\times\GHom0{G,E},E}=\Hom0{G\times\GHom0{G,E},\GHom0{G,E}}
	\]
	which is obviously a representation $R$ of $G$ on $\GHom0{G,E}$. Consider the differential $dR$. The action $dR(x,f)$ for $x\in\ger g$ and $f\in\sh O_G(G_0,E^{1|1})$ is denoted by $xf$. Its extension to $\Uenv0{\ger g}$ is denoted similarly. 

	We now define a morphism $\vphi:G_0\times\ger g_1\to G$: Let $\vphi_0=\id_{G_0}$, and define $\vphi^*$ by 
	\[
		\vphi^*(f)(u;g)=(-1)^{\Abs0u\Abs0f}(uf)(g)
	\]
	for $f\in\sh O_G(V, E_U)$, $u\in\Uenv0{\ger g}$, and $g\in V$. 
\end{Par}

The following proposition is well-known, \cf \cite{a-hc,vishnyakova}.

\begin{Prop}[sgrp-pair-eq]
	For any finite-dimensional Lie supergroup, $\vphi:G\to G_0\times\ger g_1$ is an isomorphism of supergroups. The categories of finite-dimensional Lie supergroups and finite-dimensional supergroup pairs are equivalent. 
\end{Prop}

We also have the following result, which is well-known for the case of finite-dimensional representations.

\begin{Prop}[repn-sgrp-pair]
	Let $G$ be a finite-dimensional Lie supergroup with associated supergroup pair $(\ger g,G_0)$. The categories of $G$-representations and smooth $(\ger g,G_0)$-representations are isomorphic. 
\end{Prop}

\begin{proof}
	Clearly, a representation of $G$ defines a representation of $(\ger g,G_0)$, and similarly for morphisms. Conversely, for a representation $(d\alpha,\alpha_0)$ of $(\ger g,G_0)$, define $\alpha:G\times V\to V$ to be the element of 
	\[
		\Hom0{G,\GEnd0V}=\sh O_G(G_0,\GEnd0V)\subset\Hom[_{\Uenv0{\ger g_0}}]0{\Uenv0{\ger g},\Ct[^\infty]0{G_0,\GEnd0V}}
	\]
	given by 
	\[
		\alpha(u;g)(v)=(\alpha_0(g)\circ d\alpha(u))(v)\ .
	\]
	One easily checks that this is well-defined, and inverse to the first construction. 
	
	If $f:V\to W$ is a continuous even linear map which is $(\ger g,G_0)$-equivariant, where $(W,\beta)$ is another $(\ger g,G_0)$-representation, then $f_0$ is $G_0$-equivariant, and for $u\in F^p\Uenv0{\ger g}$, $p\sge1$, one has by the definition of $\circ$
	\begin{align*}
		(f\circ\alpha)(u\otimes 1;g,v)&=f\Parens1{\alpha(u\otimes1;g,v);\alpha_0(g)(v)}\\
		&=f\Parens1{(\alpha_0(g)\circ d\alpha(u))(v)}\\
		&=\Parens1{\beta_0(v)\circ d\beta(u)}(f(v))=(\beta\circ f)(u\otimes 1;g,v)
	\end{align*}
	because $\alpha$ and $\beta$ are linear over $G$. Similarly, for $w\in V$, 
	\begin{align*}
		(f\circ\alpha)(u\otimes w;g,v)&=f\Parens1{\alpha(u\otimes w;g,v);\alpha_0(g)(v)}\\
		&=f\Parens1{(\alpha_0(g)\circ d\alpha(u))(w)}\\
		&=\Parens1{\beta_0(v)\circ d\beta(u)}(f(w))=(\beta\circ f)(u\otimes w;g,v)\ .
	\end{align*}
	Finally, for $p\in S^k(V)$, $k>1$, we have
	\[
		(f\circ\alpha)(u\otimes p;g,v)=0=(\beta\circ f)(u\otimes p;g,v)\ .
	\]
	This shows that $f$ is a morphism of $G$-representations. 
\end{proof}

\section{Applications}

\subsection{Loop supergroups and superloop supergroups}

The purpose of this section is to give a global model for a supergroup of `loops' $\mathbb S^1\to G$ or `superloops' $\mathbb S^{1|1}\to G$. We do not go into the subject in any depth, and certainly our considerations can only be considered as a very first step towards an understanding of such objects. However, we feel that even these basic results illustrate the utility of the category $\SMan$. 

In what follows, we fix a finite-dimensional Lie supergroup $G$. (The finite-dimensionality could be replaced by local metrisability.)

\begin{Def}
	We let $LG=\GHom0{\mathbb S^1,G}$ and $L^{1|1}G=\GHom0{\mathbb S^{1|1},G}$ where the supercircle $\mathbb S^{1|1}=\mathbb S^1\times\reals^{0|1}$. We call $LG$ the \Define{supergroup of loops in $G$} and $L^{1|1}$ the \Define{supergroup of superloops in $G$}. 
	
	Similarly, we set $L\ger g=\GHom0{\mathbb S^1,\ger g}$ and $L^{1|1}\ger g=\GHom0{\mathbb S^{1|1},\ger g}$, and call these the \Define{loop superalgebra} and the \Define{superloop superalgebra}, respectively. 
\end{Def}

\begin{Prop}
	The supermanifolds $LG$ and $L^{1|1}G$ are naturally Lie supergroups, and $L^{1|1}G$ is naturally isomorphic, as a Lie supergroup, to $L(\Pi TG)$. Here, $\Pi TX=\GHom0{\reals^{0|1},X}=T^AX$ for $A=\reals[\theta]$, $\theta$ being odd.  The Lie superalgebra of $LG$ is $L\ger g$, and the Lie superalgebra of $L^{1|1}G$ is $L^{1|1}\ger g$. 
\end{Prop}

\begin{proof}
	As a consequence of their definition, inner homs $\GHom0{X,-}$ commute with products:
	\begin{align*}
		\Hom0{S,\GHom0{X,Y\times Z}}&=\Hom0{S\times X,Y\times Z}\\
		&=\Hom0{S\times X,Y}\times\Hom0{S\times X,Z}\\
		&=\Hom0{S,\GHom0{X,Y}}\times\Hom0{S,\GHom0{X,Z}}\ .
	\end{align*}
	Hence, $LG$ and $L^{1|1}G$ are Lie supergroups with structure morphisms $L(m)$, \etc
	
	The isomorphism $L^{1|1}G\to L(\Pi TG)$ is a consequence of 
	\begin{align*}
		\Hom0{S,L^{1|1}G}&=\Hom0{S\times\mathbb S^1\times\reals^{0|1},G}=\Hom0{S\times\mathbb S^1\times\Spec\reals[\theta],G}\\
		&=\Hom0{S\times\mathbb S^1,\Pi TG}=\Hom0{S,L(\Pi TG)}\ ,
	\end{align*}
	together with the fact that the Lie supergroup structure on $\Pi TG$ is given by applying $\Pi T(-)$ to the structure morphisms of $G$. 
	
	Next, we compute $T(LG)$ and $T^{1|1}(LG)$. To that end, observe that inner homs commute whenever they are defined, 
	\begin{align*}
		\Hom0{S,\GHom0{X,\GHom0{Y,Z}}}&=\Hom0{S\times X,\GHom0{Y,Z}}\\
		&=\Hom0{S\times X\times Y,Z}=\Hom0{S\times Y\times X,Z}\\
		&=\Hom0{S,\GHom0{Y,\GHom0{X,Z}}}\ .
	\end{align*}
	In particular, this shows that $T\GHom0{X,Y}=\GHom0{X,TY}$. A variant of the above considerations shows that inner homs commute with fibred products, so that $T(LG)=L\ger g$ and $TL^{1|1}G=L^{1|1}\ger g$. 
\end{proof}

\begin{Par}
 	We can determine the supergroup pair associated to $LG$ explicitly, as follows: 
	\[
		(LG)_0=\Hom0{\mathbb S^1,G}=\Ct[^\infty]0{\mathbb S^1,G_0}\ ,\ L\ger g=\GHom0{\mathbb S^1,\ger g}=\Ct[^\infty]0{\mathbb S^1,\ger g}
	\]
	with their usual Fr\'echet topology. For the latter, the grading is induced from $\ger g$. 
	
	In the case of the superloops, we have $G\cong G_0\times\ger g_1$ by \thmref{Prop}{sgrp-pair-eq}, so $\Pi TG=G\times\Pi\ger g=G_0\times\ger g_1\times\Pi\ger g_0\times\Pi\ger g_1$. Thus, 
	\[
		(L^{1|1}G)_0=\Ct[^\infty]0{\mathbb S^1,G_0\times\Pi\ger g_1}\ ,\ L^{1|1}\ger g=\Ct[^\infty]0{\mathbb S^1,\ger g\oplus\Pi\ger g}\ .
	\]
	
	This supergroup pair is given by applying $\Ct[^\infty]0{\mathbb S^1,-}$ to the supergroup pair $(G_0\times\Pi\ger g_1,\ger g\oplus\Pi\ger g)$. Here, the adjoint action of $G_0\times\ger g_1$ on $\ger g\oplus\Pi\ger g$ is given by
	\[
		\Ad((g,x),y\oplus z)=\Ad(g)(y)\oplus\Parens1{\Ad(g)([x,y])+\Ad(g)(z-x)}\ ,
	\]
	and the bracket of $\ger g\oplus\Pi\ger g$ is
	\[
		[x_1\oplus y_1,x_2\oplus y_2]=[x_1,x_2]\oplus\Parens1{[x_1,y_2]+[y_1,x_2]}\ .
	\]
	All of the latter can be derived in the usual category of super-ringed spaces, by considering $S\times\Spec\reals[\theta]$-points, so we don't give the details of the derivation. 
\end{Par}

\subsection{Induced representations}

In this section, we give some elements of induced representations. As before, these are only the basics, and they serve as an illustration of the category $\SMan$. In what follows, all Lie supergroups will be finite dimensional.

\begin{Def}
	Let $P/S$ be a bundle with fibre $H$, where $H$ is a Lie supergroup. Assume there is a right $H$-action on $P$ \scth the bundle projection is equivariant (for the trivial action on $S$), and that there is an atlas of $H$-equivariant trivialisations. Then $P/S$ is called a \Define{principal $H$-bundle}.
\end{Def}

\begin{Lem}
	Let $P/S$ be a principal $H$-bundle and $V$ an $H$-representation. Then the coequaliser $P\times^HV$ of $\xymatrix@1{X\times H\times V\ar@<+0.5ex>[r]\ar@<-0.5ex>[r]&X\times V}$ exists, and is a vector bundle on $S$ with fibre $V$. In particular, if $H$ is a closed subsupergroup of $G$, then $G\times^HV$ exists and is a $G$-equivariant vector bundle over the quotient $G/H$. 
\end{Lem}

\begin{proof}
	The first statement is easily proved for a trivial bundle, and the general case follows by patching bundles. As for the latter statement, it is known from Ref.~\cite{kostant} that $G/H$ exists and $G\to G/H$ is a principal $H$-bundle. The equivariance is easy to check. 
\end{proof}

\begin{Prop}[ind-def]
	Let $H$ be a closed subsupergroup of $G$ and $V$ a metrisable $H$-representation. Then 
	\[
		\ind_H^G(V):=\Gamma(G/H,G\times^HV)
	\]
	carries a natural $G$-representation.
\end{Prop}

\begin{proof}
	For any $S\in\Ob\SMan$, $G/H(S)=\Hom[_S]0{S,S\times G/H}$ and 
	\[
		(\ind_H^GV)(S)=\Hom[_{G/H}]0{S\times G/H,G\times^HV}\ ,
	\]
	so for $x\in_SG/H$, $f\in_S\ind_H^G(V)$, we have $f(x)\in_SG\times^HV$. Thus, for $g\in_SG$, 
	\[
		(g.f)(x)=g.(f(g^{-1}.x))\in_SG\times^HV\ ,
	\]
	and for $\pi=\pi_{G\times^HV}$, 
	\[
		\pi\Parens1{g(f(g^{-1}.x))}=g.\pi(f(g^{-1}.x))=g.(g^{-1}.x)=x\ ,
	\]
	so this indeed defines $g.f\in_S\ind_H^GV$, and a morphism $G\times\ind_H^GV\to\ind_H^GV$. Standard computations show that it is indeed a representation. 
\end{proof}

\begin{Prop}
	Let $H_1\to H_2\to G$ be a chain of closed subsupergroups and $V$ a metrisable $H_1$-representation. Then there is a natural isomorphism 
	\[
		\ind_{H_1}^GV\cong\ind_{H_2}^G\ind_{H_1}^{H_2}V
	\]
	of $G$-representations.
\end{Prop}

\begin{proof}
	Clearly, there is an isomorphism of $G$-equivariant vector bundles
	\[
		G\times^{H_1}V\cong G\times^{H_2}(H_2\times^{H_1}V)\ .
	\]
	By applying the definition of the section functor, one easily obtains the claim. 
\end{proof}

We have the following version of Frobenius reciprocity. 

\begin{Prop}[frobenius]
	Let $h$ be a closed subsupergroup of $G$, $V$ a metrisable $H$-representa\-tion and $W$ a $G$-representa\-tion. There is a natural linear isomorphism
	\[
		\Hom[^G]0{W,\ind_H^GV}\cong\Hom[^H]0{W|_H,V}\ .
	\]
	Here, $\Hom[^G]0{\cdot,\cdot}$ denotes the space of $G$-equivariant even linear morphisms, and similarly for $H$. 
\end{Prop}

\begin{proof}
	We fix the standard base point $o\in(G/H)_0=G_0/H_0$. Then $V$ identifies with the fibre of $G\times^HV$ at $o$, and as in the proof of \thmref{Prop}{ind-def}, evaluation at $o$ gives a linear morphism $\ind_H^GV\to V$. We thus obtain a linear map
	\[
		\vphi:\Hom[^G]0{W,\ind_H^GV}\to\Hom[^H]0{W|_H,V}
	\]
	by $\vphi(f)(w)=f(w)(o)$. Indeed, for $h\in_SH$, $w\in_SW$, 
	\begin{align*}
		\vphi(f)(h.w)&=f(h.w)(o)=(h.f(w))(o)\\
		&=h.(f(w)(h^{-1}.o))=h.(f(w)(o))=h.\vphi(f)(w)\ ,
	\end{align*}
	so $\vphi$ is well-defined. 
	
	Conversely, we define a map $\psi$. Let $f\in\Hom[^H]0{W|_H,V}$ and $w\in_S W$. Define $k\in_S\ind_H^GV$ by 
	\[
		k(g.o)=g.(f(g^{-1}.w))\in_SG\times^HV\ , 
	\]
	where the canonical map $V=(G\times^HV)_o\to G\times^HV$ is suppressed in the notation. The stabliser of $G(S)$ at $o$ is $H(S)$, and a standard check shows that this is a well-defined $S$-point of $G\times^HV$. Then set $\psi(f)(w)=k$. This is well-defined, since
	\begin{align*}
		\psi(f)(g_1.w)(g_2.o)&=g_2.(f(g_2^{-1}g_1.w))=g_1.(g_1^{-1}g_2.(f((g_1^{-1}g_2)^{-1}.w)))\\
		&=g_1.(\psi(f)(w)(g_1^{-2}g_2.o))=(g_1.\psi(f)(w))(g_2.o)
	\end{align*}
	\fa $g_1,g_2\in_SG$, $w\in_SW$. The linearity and $H$-equivariance of $\psi$ are no harder to check. Moreover, it is clear that $\vphi$ and $\psi$ are mutually inverse. 
\end{proof}

\begin{Rem}
	Induced representations can be considered in terms of supergroup pairs, \cf \cite{chemla} for the finite-dimensional case. The above results have been derived in these terms and applied in the proof of a super version of the Cartan--Helgason Theorem, \vq \cite{schmittner-dip}. 
\end{Rem}

\subsection{Convolution algebras}

As above, we will assume all Lie supergroups $G$ to be finite-dimensional; moreover, $G_0$ should be metrisable. 

\begin{Def}
	For $X\in\SMan$ locally in $\cat{fd}$, and $E$ a locally convex super-vector space (over $\reals$ or $\cplxs$), define $\sh E'(X,E')$ to be the topological dual of the locally convex super-vector space $\GHom0{X,E}=\sh O_X(X_0,E^{1|1})$, endowed with the strong topological and the induced grading. In particular, for $E=\cplxs$, we write $\sh E'(X)=\sh E'(X,\cplxs)$. We call the members of $\sh E'(X,E')$ \Define{compactly supported $E$-distributions} and those of $\sh E'(X)$ \Define{compactly supported distributions}. 
	
	An important subset is formed by $\Ct[_c^\infty]0G=\Gamma_c(G,\ABer0{TG})\otimes\cplxs$, the compactly supported sections of the Berezinian density bundle, \cf \cite{leites,ah-int}. These are called \Define{smooth densities} (with compact support), and embed into $\sh E'(X)$ \via
	\[
		\Dual0f\omega=\int_Xf\cdot\omega\mathfa f\in\sh O_X(X_0,\cplxs^{1|1})\,,\,\omega\in\Ct[_c^\infty]0X\ .
	\]
	
	Let $G$ be a Lie supergroup and $\alpha:G\times X\to X$ an action. For $\mu\in\sh E'(G)$, $\nu\in\sh E'(X,E')$, define $\mu*\nu\in\sh E'(X,E')$ by 
	\[
		\Dual0f{\mu*\nu}=\Dual0{\alpha^*f}{\mu\otimes\nu}
	\]
	\fa $f\in\GHom0{X,E}=\sh O_X(X_0,E^{1|1})$. 
	
	Let $V$ be a locally convex super-vector space and $\alpha:G\times V\to V$ a representation. For $\mu\in\sh E'(G)$ and $v'\in V'$, define $\alpha(\mu)v'\in V'$ by
	\[
		\Dual0v{\alpha(\mu)v'}=\Dual0{\alpha^*(v)}{\mu\otimes v'}
	\]
	\fa $v\in V$. Here, observe that $V\subset\sh O_V(V_0,\cplxs^{1|1})$ and $\alpha^*(V)\subset\sh O_G(G_0,\cplxs^{1|1})\otimes V$ because $\alpha$ is linear over $G$. 
\end{Def}

\begin{Prop}
	With $*$ as a product, $\sh E'(G)$ is a locally $m$-convex associative unital superalgebra over $\cplxs$. The operations $*$ on $\sh E'(X,E')$ and $\alpha(-)$ on $V'$ define unital topological $\sh E'(G)$-supermodule structures.
\end{Prop}

\begin{proof}
	Using partitions of unity and \thmref{Lem}{smooth-curry}, it follows that the locally convex super-vector space $\sh O_{G\times X}(G_0\times X_0,E^{1|1})$ contains $\sh O_G(G_0,\cplxs^{1|1})\otimes\sh O_X(X_0,E^{1|1})$ as a dense subspace. Therefore, $\mu*\nu$ is well-defined for $\mu\in\sh E'(G)$ and $\nu\in\sh E'(X,E')$. The remaining statements are straightforward to check. 
\end{proof}

\begin{Def}
	Let $G$ be a Lie supergroup. Let $V$ be a complete locally convex super-vector space and $\alpha:G\times V\to V$ a representation. 
	
	Let $V'$ be the topological dual of $V$, endowed with the strong topology (of bounded convergence on bounded subsets). Endowed with the action given by 
	\[
		\Dual0v{\check\alpha_0(g)v'}=\Dual0{\alpha_0(g^{-1})v}{v'}\nd\Dual0v{d\check\alpha(x)v'}=-(-1)^{\Abs0x\Abs0v}\Dual0{d\alpha(x)v}{v'}\ ,
	\] 
	it is a smooth representation of $(\ger g,G_0)$, as follows directly from Ref.~\cite[Propositions 2.3-4, and discussion following Definition 2.1]{bruhat-thesis}. Thus, $V'$ is naturally a $G$-representation. 
	
	Let $V$ be semi-reflexive and complete, and $\alpha_0:G\times V\to V$ a continuous representation of $G_0$. Let $V_\infty$ be the set of smooth vectors with its natural topology. It is a Fr\'echet space, and by \loccit, $V'_{-\infty}=(V_\infty)'$ is a smooth $G_0$-representation containing the continuous $G_0$-representation $V'$. If $V$ is, moreover, reflexive, we may define $V'_\infty$ and $V_{-\infty}=(V'_\infty)'$, which, by the density of $V'_\infty$ in $V'$, contains $V$ as a subspace. 
	
	If $V$ is a continuous $G_0$-representation and the representation of $(G_0,\ger g_0)$ on $V_\infty$ extends to a smooth $(\ger g,G_0)$-representation, then we say that $V$ is a \Define{continuous $(\ger g,G)$-representation}. This is an immediate generalisation of the definitions given in Refs.~\cite{cctv,merigon-neeb-salmasian} for the case of unitary representations. 
\end{Def}

We have the following generalisation of the Theorem of Dixmier--Malliavin. 

\begin{Th}[super-dixmal]
	Let $G$ be a Lie supergroup and $\alpha_0$ a $G_0$-representation on a reflexive Fr\'echet super-vector space $V$. Assume that the $(\ger g_0,G_0)$-representation on $V_\infty$ extends to a representation $G$ of $V_\infty$. Then $\Ct[^\infty_c]0G\subset\sh E'(G)$ acts on $V$, and 
	\[
		\alpha(\Ct[^\infty_c]0G)V=\alpha(\Ct[^\infty_c]0G)V_\infty=V_\infty\ .
	\]
	
	More precisely, given $v\in V_\infty$, there exist finitely many $f_j\in\Ct[^\infty_c]0G$ and $v_j$ in the closed $(G_0,\ger g)$-invariant subspace of $V_\infty$ generated by $v$, \scth $v=\sum_j\alpha(f_j)v_j$. 
\end{Th}

In the \emph{proof}, we need the following lemma.

\begin{Lem}
	Let $\alpha_0$ be a continuous $G_0$-representation on a reflexive Fr\'echet space $V$. If $f\in\Ct[_c^\infty]0{G_0}$, then $\alpha_0(f)$ is a continuous linear map $V_{-\infty}\to V_\infty$. 
\end{Lem}

\begin{proof}
	It is known that $\check\alpha_0(V')\subset V'_\infty$. For $v\in V_{-\infty}$, define $\alpha_0(f)v\in V$ by 
	\[
		\Dual0{v'}{\alpha_0(f)v}=\Dual0{\check\alpha_0(i^*f)}v
	\]
	where $i(g)=g^{-1}$. Since $i$ is a diffeomorphism and $V$ is semi-reflexive, it is clear that $\alpha_0(f)v$ is well-defined and continuous in $v$, and that this definition extends $\alpha_0(f)$ from $V$ to $V_{-\infty}$. 
	
	By Dixmier--Malliavin \cite{dixmier-malliavin}, there exist finitely many $f_j,h_j\in\Ct[^\infty_c]0{G_0}$ \scth $f=\sum_jf_j*h_j$. Then 
	\[
		\alpha_0(f)v=\sum\nolimits_j\alpha_0(f_j)\alpha_0(h_j)v
	\]
	Since $\alpha_0(h_j)v\in V$ and $\alpha_j(f_j)(V)\subset V_\infty$, we conclude that $\alpha_0(f)(V_{-\infty})\subset V_\infty$. 
	
	To see the continuity, let $p$ be a continuous seminorm on $V$ and $u\in\Uenv0{\ger g_0}$. Then 
	\[
		p(d\alpha_0(u)\alpha_0(f)v)=p(\alpha_0(uf)v)
	\]
	depends continuously on $v\in V_{-\infty}$. Since the seminorms $p\circ d\alpha_0(u)$ define the topology on $V_\infty$, this proves the claim. 
\end{proof}

\begin{proof}[\protect{Proof of \thmref{Th}{super-dixmal}}]
	We observe that $G$ acts on $V_{-\infty}$, so $\alpha(\Ct[^\infty_c]0G)V\subset V_{-\infty}$. 
	
	Let $q=\dim\ger g_1$. There exist, for $I\subset\{1,\dotsc,q\}$, $u_I,v_I\in F^q\Uenv0{\ger g}$ \scth \fa $f\in\Ct[^\infty_c]0G$, $v\in V$, we have 
	\[
		\alpha(f)v=\sum\nolimits_I\alpha_0(u_If)d\alpha(v_I)v\in\alpha_0(\Ct[^\infty_c]0{G_0})V_{-\infty}\subset V_\infty\ .
	\]
	In particular, $\Ct[^\infty_c]0G$ acts on $V$, and maps $V$ to $V_\infty$. 
	
	To prove the theorem, we need to show that $V_\infty\subset\alpha(\Ct[_c^\infty]0G)V_\infty$. To that end, we introduce the following terminology: A closed Lie subsupergroup $H$ is called \Define{singly generated} if its Lie superalgebra $\ger h$ is generated as a Lie superalgebra by a single homogeneous element. Then the following is clear: Any singly generated Lie subsupergroup is locally isomorphic to $\reals$, to $\reals^{0|1}$, or to $\reals^{1|1}$ where the Lie superalgebra has the unique non-zero homogeneous relation $x=[y,y]$. Moreover, there exist singly generated closed Lie subsupergroups $H_1,\dotsc,H_n$ \scth the multiplication morphism $m:H_1\times\dotsm\times H_n\to G$ is an isomorphism in a neighbourhood $U$ of the identity. 
	
	Now, fix $v\in V_\infty$. We claim the following: For any singly generated closed Lie subsupergroup $H$, there exist $f_1,f_2\in\Ct[_c^\infty]0{H}\subset\sh E'(H)\subset\sh E'(G)$ and $w\in V_\infty$ with $\supp f_j\subset H\cap U$, \scth $v=\pi(f_0)v+\pi(f_1)w$. Since $\reals\subset\reals^{1|1}$ as a closed Lie subsupergroup, this follows from Dixmier--Malliavin \cite[]{dixmier-malliavin} in case $H_0$ is locally isomorphic to $\reals$. 
	
	In case $H$ is isomorphic to $\reals^{0|1}$, we have $\sh O_H(H_0,\reals^{1|1})=\reals[\tau]$ where $\tau$ is odd. It follows that $\int_H f\,D\tau=\frac d{d\tau}f$ is a smooth density, and $\int_Hf\tau\,D\tau=f(0)$. Thus, $\delta=\tau\,D\tau$ is smooth, and hence, the statement is obvious in this case. 
	
	Applying the statement inductively, we find $f_0^j,f_1^j\in\Ct[^\infty_c]0{H_j}\subset\sh E'(G)$ and $w_{i_1,\dotsc,i_n}\in V_\infty$, $i_j=0,1$, \scth 
	\[
		v=\sum_{i_1,\dotsc,i_n=0,1}\alpha(f^1_{i_1}*\dotsm*f_{i_n}^n)w_{i_1,\dotsc,i_n}\ .
	\]
	Now, for $f_j\in\Ct[^\infty_c]0{H_j}$, $f\in\sh O_G(G_0,\cplxs^{1|1})$, we have 
	\[
		\Dual0f{f_1*\dotsm*f_n}=\Dual0{m^*f}{f_1\otimes\dotsm\otimes f_n}
	\]
	where $m:H_1\times\dotsm\times H_n\to G$ is the $n$-fold multiplication morphism. Since $f_1\otimes\dotsm\otimes f_n$ is smooth, so is $f_1*\dotsm*f_n$. This proves the claim. 
\end{proof}

\begin{Rem}
	In the proof of their theorem, Dixmier--Malliavin use the fact that the exponential map is a group homomorphism for one-parameter Lie groups. The corresponding fact for Lie supergroups is false. The singly generated closed Lie subsupergroups used above are a replacement for the one-parameter groups. 
\end{Rem}

\begin{Lem}[approx-unit]
	Let $G$ be a Lie supergroup and $\sh U$ a basis of relatively compact neighbourhoods of $1$ in $G_0$. For every $U\in\sh U$, there exists $\psi_U\in\Ct[^\infty_c]0G$ \scth $\supp\psi_U\subset U$, $i^*\psi_U=\psi_U$, $\int_G\psi_U=1$, and $\psi_U\to\delta$ in $\sh E'(G)$. 
\end{Lem}

\begin{proof}
	The first statment is local, and it is not hard to construct $\psi_U$ \scth $\int_G\psi_U=1$. Since $i:G\to G$ is an isomorphism, $\int_Gi^*\psi_U=1$, and we may replace $\psi_U$ by $\frac12(\psi_U+i^*\psi_U)$. As to the second statement, 
	\[
		\int_Gf\psi_U-f(1)=\int_G(f-f(1))\psi_U\ ,
	\]
	so to prove convergence, one may assume that $\supp f$ is compact and $f(1)=0$. Since $f$ is compactly supported, $(uf)|_{G_0}$ is left uniformly continuous \fa $u\in\Uenv0{\ger g}$. 
	
	Moreover, with $q=\dim\ger g_1$, there exists a basis $\sh B$ of $F^q\Uenv0{\ger g}$ \scth $1\in\sh B$, $\eps(u)=0$ for $u\in\sh B$, $1\neq u$, and 
	\[
		\Abs2{\int_Gh\omega}\sle\sum_{u_1,u_2\in\sh B}\Abs2{\int_G u_1\omega}\cdot\sup_{g\in G_0}\Abs0{u_2h(g)}
	\]
	\fa integrable densities $\omega$ and all compactly supported superfunctions $h$. 
	
	Then for each $\eps>0$, there exists $U\in\sh U$ \scth \fa $V\in\sh U$, $V\subset U$, one has that $\Abs0{uf}(g)<\eps$ \fa $u\in\sh B$, and all $g\in V$. Since $\int_G (u\psi_V)=0$ if $\eps(u)=0$, 
	\[
		\Abs2{\int_G f\psi_V}\sle\Abs2{\int_G \psi_V}\cdot\sum_{u\in\sh B}\sup_{g\in G_0}\Abs0{uf(g)}<\#\sh B\cdot\eps\ .
	\]
	This proves the claim.
\end{proof}

\begin{Prop}[multi]
	Let $V$, $W$ be continuous reflexive Fr\'echet representations of the supergroup pair $(\ger g,G_0)$ with associated Lie supergroup $G$. 
	\begin{enumerate}
		\item There is a bijection between: 
		\begin{enumerate}[\normalfont (a).]
			\item Closed $G_0$-invariant subspaces $U\subset V$ \scth $U_\infty$ is $\ger g$-invariant; 
			\item closed $G$-invariant subspaces of $V_\infty$; and 
			\item closed $\Ct[^\infty_c]0G$-invariant subspaces of $V$. 
		\end{enumerate}
		\item There is a bijection between:
		\begin{enumerate}[\normalfont (a).]
			\item Linear $G_0$-equivariant morphisms $\phi:V\to W$ \scth the restriction $\phi:V_\infty\to W_\infty$ is $\ger g$-invariant; 
			\item linear $G$-equivariant morphisms $V_\infty\to W_\infty$; and 
			\item even continuous linear maps $V\to W$ which are $\Ct[^\infty_c]0G$-linear.
		\end{enumerate}		
	\end{enumerate}
\end{Prop}

\begin{proof}
	(i). In view of \thmref{Prop}{repn-sgrp-pair}, the bijection between (a) and (b) follows from \cite[Proposition 2.6]{bruhat-thesis}. Assume that $U\subset V$ is a closed subspace of $V$. If $U$ is $G_0$-invariant and $U_\infty$ is $\ger g$-invariant, then by \thmref{Th}{super-dixmal}, $U$ is $\Ct[^\infty_c]0G$-invariant. 
	
	Conversely, if $U$ is $\Ct[^\infty_c]0G$-invariant, then so is $U_\infty=U\cap V_\infty$. Let $u\in U$ and choose $\psi_U\in\Ct[^\infty_c]0G$ as in \thmref{Lem}{approx-unit}. Then for $g\in G_0$, 
	\[
		\alpha(g)u=\alpha(g)\alpha(\delta)u=\lim\nolimits_U\alpha(\delta_g*\psi_U)u
	\]
	where the right side is in $U$ by assumption. Similarly, if $u\in U_\infty$, then for $x\in\ger g$, 
	\[
		d\alpha(x)u=d\alpha(x)\alpha(\delta)u=\lim\nolimits_U\alpha(x\psi_U)u
	\]
	where the right hand side is $U$, and the left hand side in $V_\infty$, so that we conclude $d\alpha(x)u\in U\cap V_\infty=U_\infty$. This proves the bijection between (a) and (c).
	
	(ii). The proof of this statement is entirely analogous to that of (i), so we leave it to the reader. 
\end{proof}

\begin{Rem}
	As the above proposition shows, the study of multiplicities in $G$-representations can adressed equivalently by the study of the corresponding problem for representations of the corresponding supergroup pair, or the convolution algebra of smooth densities $\Ct[^\infty_c]0G$. 
\end{Rem}

\bibliographystyle{alpha}%
\bibliography{a-srep}%

\end{document}